\documentclass[11pt]{article}
\usepackage{amsmath}
\usepackage{amsfonts}
\usepackage{graphicx}
\usepackage{latexsym}
\usepackage{color}

\newcommand{\R}{\mathbb{R}}
\newcommand{\E}{\mathbb{E}}
\newcommand{\PP}{\mathbb{P}}

\newcommand{\III}{\mathbf{I}}

\newcommand{\ppp}{\mathbf{p}}

\newcommand{\al}{\alpha}
\newcommand{\la}{\lambda}
\newcommand{\La}{\Lambda}
\newcommand{\ga}{\gamma}
\newcommand{\Ga}{\Gamma}
\newcommand{\ka}{\kappa}

\newcommand{\si}{\sigma}

\newcommand{\te}{\theta}
\newcommand{\Te}{\Theta}
\newcommand{\be}{\beta}
\newcommand{\ep}{\varepsilon}

\newcommand{\de}{\delta}
\newcommand{\De}{\Delta}
\newcommand{\om}{\omega}
\newcommand{\Om}{\Omega}
\newcommand{\ze}{\zeta}

\newcommand{\up}{\upsilon}

\newcommand{\Up}{\Upsilon}

\newcommand{\ba}{\mathcal{B}}
\newcommand{\ca}{\mathcal{C}}

\newcommand{\f}{\mathcal{F}}
\newcommand{\g}{\mathcal{G}}

\newcommand{\ia}{\mathcal{I}}
\newcommand{\ja}{\mathcal{J}}

\newcommand{\laa}{\mathcal{L}}
\newcommand{\m}{\mathcal{M}}
\newcommand{\n}{\mathcal{N}}
\newcommand{\p}{\mathcal{P}}

\newcommand{\s}{\mathcal{S}}

\newcommand{\ua}{\mathcal{U}}

\newcommand{\y}{\mathcal{Y}}
\newcommand{\z}{\mathcal{Z}}

\newcommand{\wR}{\widehat{R}}

\newcommand{\wS}{\widehat{S}}

\newcommand{\WE}{\widetilde{\mathbb{E}}}

\newcommand{\WV}{\widetilde{V}}

\newcommand{\Wsi}{\widetilde{\sigma}}

\newcommand{\Wze}{\widetilde{\zeta}}

\newcommand{\BA}{\overline{A}{}}

\newcommand{\BF}{\overline{F}{}}

\newcommand{\BB}{\overline{B}{}}

\newcommand{\BV}{\overline{V}{}}
\newcommand{\BW}{\overline{W}{}}

\newcommand{\BZ}{\overline{Z}{}}

\newcommand{\BPhi}{\overline{\Phi}{}}

\newcommand{\BTe}{\overline{\Theta}{}}

\newcommand{\BOm}{\overline{\Om}{}}

\newcommand{\BE}{\overline{\E}{}}

\newcommand{\BUp}{\overline{\Up}{}}

\newcommand{\Bfa}{\overline{\f}{}}
\newcommand{\BPP}{\overline{\PP}{}}

\newcommand{\uu}{\underline{u}}

\def\nib{\noindent\bf}

\def\ni{\noindent}

\newcommand{\rdn}{\sqrt{\De_n}}

\newcommand{\rdnn}{\frac{1}{\rdn}}

\newcommand{\dd}{\De^n_i}

\newcommand{\vsc}{\vskip 5mm}
\newcommand{\vsd}{\vskip 2mm}
\newcommand{\vst}{\vskip 3mm}
\newcommand{\vsq}{\vskip 4mm}

\newcommand{\toop}{\stackrel{\PP}{\longrightarrow}}

\newcommand{\tolls}{\stackrel{\laa-\mbox{\tiny s}}{\Longrightarrow}}

\newcommand{\tols}{\stackrel{\laa-\mbox{\tiny s}}{\longrightarrow}}

\newcommand{\toucp}{\stackrel{\mbox{\tiny u.c.p.}}{\Longrightarrow}}

\newcommand{\sign}{\mbox{\rm sign}}

\newcommand{\proba}{(\Omega ,\f,(\f_t)_{t\geq0},\PP)}

\newcommand{\probt}{(\widetilde{\Omega},\widetilde{\f},
(\widetilde{\f}_t)_{t\geq0},\widetilde{\PP})}

\newcommand{\bee}{\begin{equation}}
\newcommand{\eee}{\end{equation}}
\newcommand{\bea}{\begin{eqnarray}}
\newcommand{\eea}{\end{eqnarray}}
\newcommand{\bean}{\begin{eqnarray*}}
\newcommand{\eean}{\end{eqnarray*}}

\newcommand{\qed}{$\hfill\Box$}

\newtheorem{prop}{Proposition}[section]
\newtheorem{cor}[prop]{Corollary}

\newtheorem{lem}[prop]{Lemma}

\newtheorem{theo}[prop]{Theorem}
\newtheorem{rem}[prop]{Remark}

\begin{document}

\title{A test for the rank of the volatility process: the random perturbation
approach}
\author{Jean Jacod \thanks{Institut de Math\'ematiques de Jussieu,
4 Place Jussieu, 75 005 Paris, France (CNRS -- UMR 7586, and
Universit\'e Pierre et
Marie Curie), Email: jean.jacod@upmc.fr} \and
Mark Podolskij \thanks{Department of Mathematics, Heidelberg University,
INF 294, 69120 Heidelberg,
Germany, Email: m.podolskij@uni-heidelberg.de.}}

\date{\today}

\maketitle

\begin{abstract}
In this paper we present a test for the maximal rank of the matrix-valued
volatility process in the continuous
It\^o semimartingale framework. Our idea is based upon a random perturbation
of the original high frequency observations
of an It\^o semimartingale, which opens the way for rank testing. We develop
the complete limit theory for the test statistic
and apply it to various null and alternative hypotheses. Finally, we
demonstrate a homoscedasticity test for the rank process.

\ \

\bigskip


\end{abstract}

\section{Introduction}
\label{Intro}
\setcounter{equation}{0}
\renewcommand{\theequation}{\thesection.\arabic{equation}}
In the last years asymptotic statistics for high frequency observations has
received a lot of attention in the literature. This interest was mainly
motivated by  financial applications, where observations of stocks or
currencies are available at very high frequencies. As under the no-arbitrage
condition prices processes must be semimartingales (see e.g. \cite{DeS}), a
lot of research has been devoted to statistics of high frequency data of
semimartingales. We refer to a recent book \cite{JP12} for a comprehensive
study of infill asymptotic for semimartingales.

This paper is devoted to testing for the maximal rank of the matrix-valued
volatility process in the continuous It\^o semimartingale framework, and
more specifically for a $d$-dimensional continuous It\^o semimartingale
$X$ which is observed at equidistant
times over a fixed time interval $[0,T]$: we observe
$(X_{i\De_n})_{0\leq i\leq [T/\De_n]}$, and the high-frequency approach
consists in assuming $\De_n\rightarrow 0$.

A continuous It\^o semimartingale can be written as
\bee\label{0-1}
dX_t = b_t dt + \si_t dW_t,
\eee
where $W$ is a Brownian motion, and there are many representations of
this form, with different Brownian motions $W$ and, accordingly, different
volatility processes $\si$. What is ``intrinsic'' is the drift coefficient
$b_t$ and the diffusion coefficient (``squared volatility'') $c_t=
\si_t\si_t^*$, in the sense that they are uniquely determined by $X$, up
to a Lebesgue-null set of times (throughout the paper $\si_t^*$ denotes the transpose of the matrix $\si$).

For modeling purposes and economical interpretation we would like to find,
and often choose, the smallest possible dimension of the Brownian motion $W$
in the representation \eqref{0-1}. Assuming further that $t\mapsto c_t$
is continuous, this smallest possible dimension is the
supremum in time of the rank of the $\R^{d \times d}$-valued process $c$
over the time interval $[0,T]$. We are further interested in homoscedasticity
testing for the rank process.

A partial answer to this question was given in \cite{JLT}. The authors of this paper studied the problem of testing 
the null hypothesis $\sup_{t\in [0,T)} \text{rank} (c_t) \geq r_0$ against $\sup_{t\in [0,T)} \text{rank} (c_t) < r_0$
for a given number $r_0$. However,
their method  does not extend to testing null hypotheses of other types, e.g. $\sup_{t\in [0,T)} \text{rank} (c_t) = r_0$ against $\sup_{t\in [0,T)} \text{rank} (c_t) \not = r_0$ (which is much more useful). In the classical setting of i.i.d or weakly dependent data various estimation methods for the rank of an unknown covariance matrix (and related objects) have been proposed. We would
like to mention Gaussian elimination method with complete pivoting of \cite{CD} and the test suggested in \cite{RS} among others.
Unfortunately, these procedures can not be applied to our statistical problem as the probabilistic structure of the process $X$
is more complex and the rank is time-varying.

Our method is based upon a random perturbation of the original data and determinant expansions. The main idea can be described as
follows: if we compute $\det(c_t + he_t)$ for a positive definite $d\times d$ matrix $e_t$ independent of $c_t$ and $h\downarrow 0$,
then, under appropriate conditions, its rate of decay to $0$ depends on the unknown rank of $c_t$. Hence, the ratio $\det(c_t + 2he_t)/\det(c_t + he_t)$
asymptotically identifies the rank of $c_t$. Indeed, our main statistic is a partial sum of squared determinants of matrices
build from $d$ consecutive increments of the process $X$ and the random perturbation is performed by a properly scaled Brownian motion $W'$, which is independent of all ingredients of $X$. We remark that perturbation methods (and matrix expansions as well)
find applications in various fields of mathematics; we refer for instance to \cite{KN} whose authors apply
matrix perturbation methods to determine the number of components in a linear mixture model from high dimensional noisy samples.
Furthermore, the methods of \cite{CM} also rely upon a generation of a new Brownian motion $W'$ although
in a completely different setting.

The paper is structured as follows. Section \ref{SAN} is devoted to  model
assumptions, testing hypotheses and test statistics.
We present the asymptotic theory for our estimators and apply it to maximal rank testing in section \ref{Res}.
In section \ref{sec4} we develop a test for the null hypothesis of constant rank. All proofs
are deferred to section \ref{P}.

\section{Model, assumptions and a random perturbation}
\label{SAN}
\setcounter{equation}{0}
\renewcommand{\theequation}{\thesection.\arabic{equation}}

\subsection{The setting and testing hypotheses} \label{sec2-1}
Our process of interest is a $d$-dimensional continuous It\^o
semimartingale $X$, given on some filtered probability space $\proba$.
In vector form, and with $W$ denoting a $q$-dimensional Brownian motion,
it can be written as
\bee\label{1-0}
X_t= X_0+ \int_0^t b_s \,ds + \int_0^t \si_s\, dW_s,
\eee
where $b_t$ is a $d$-dimensional drift process and $\si_t$ is a
$\R^{d \times q}$-valued volatility process, assumed to be continuous
in time (and indeed much more, see Assumption (H) below). We set
\bee\label{1-5}
c_t=\si_t\si_t^*,\qquad r_t=\text{rank}(c_t),\qquad
R_t=\sup_{s\in[0,t)}\,r_s.
\eee

We remark that the maximal rank $R_T$ is not bigger than the rank of the
integrated volatility $\int_0^T c_t dt$, but may be strictly smaller. As
already mentioned, it is suitable to use the smallest possible
dimension for $W$, on the time interval $[0,T]$. This is the $\PP$-essential
supremum of $\om\mapsto R_T(\om)$, but, since a single path
$t\mapsto X_t(\om)$ is (partially) observed, the only available information is
$R_T$ itself. So the problem really boils down to finding the behavior of
the process $r_t$, and for this the choice of the dimension of $W$ in
(\ref{1-0}) is irrelevant.

The rank $r_t$ is the biggest integer $r\leq d$ such that the sum of the
determinants of the matrices $(c_t^{ij})_{i,j\in J}$, where $J$ runs through
all subsets of $\{1,\cdots,d\}$ with $r$ points, is positive (with the
convention that a $0\times0$ matrix has determinant $1$); see e.g. \cite[Lemma 3]{JLT}. Since $c_t$ is
continuous, this implies that for any $r$ the random set
$\{t: r_t(\om)> r\}$ is open in $[0,T)$, so the mapping $t\mapsto r_t$ is lower
semi-continuous. In particular, the set $\{t\in[0,T):r_t(\om)=R_T(\om)\}$ is
a non-empty open subset. These properties also yield that the process $r_t$ is
predictable and that the following subsets of $\Om$, which later will
be the ``testing hypotheses'',  are $\f_T$-measurable:
\bee\label{1-2}
\begin{array}{l}
\Om^r_T~=~\{\om:\,R_T(\om)=r\}\\[1.5 ex]
\Om_T^{=}~=~\{\om:\,r_t(\om)=R_T(\om)~\text{for all}~ t\in[0,T]\}\\ [1.5 ex]
\Om^{\neq}_T~=~\{\om:~t\mapsto r_t(\om)~\text{has finitely many discontinuities
and is}\\\hskip5cm\text{ not Lebesgue-a.s. constant on $[0,T]$}\}.
\end{array}
\eee
Notice that we impose that $r_T=R_T$ in $\Om^=_T$, whereas the
lower semi-continuity only implies in general that $r_T\leq R_T$.
Observe also that {\em a priori} $t\mapsto r_t$ may be Lebesgue-a.s. constant
and still have discontinuities (even infinitely many) on $[0,T]$. So,
$\Om_T^=$ and $\Om^{\neq}_T$ are disjoint but $\Om_T^=\cup\Om^{\neq}_T
\neq\Om$ in general. The main aim of this paper is testing the null
hypothesis $\Om^r_T$ against $\Om_T^{\neq r}=\cup_{r'\neq r,0\leq r'\leq d}
\Om_T^{r'}$ (and related hypotheses) and testing the null hypothesis of
$\Om_T^{=}$ against $\Om^{\neq}_T$.

\subsection{Matrix perturbation} \label{sec2-2}
In order to explain the main idea of our method, we need to introduce
some notation. Recall that $d$
and $q$ are the dimensions of $X$ and $W$, respectively. Then $\m$
is the set of all $d\times d$ matrices, $\m_r$ for
$r\in\{0,\cdots,d\}$ is the set of all matrices in $\m$ with rank $r$,
and $\m'$ is the set of all $d\times q$ matrices. For any matrix $A$
we denote by $A_i$ the $i$th column of $A$; for any vectors
$x_1,\ldots,x_d$  in $\R^d$, we write mat$(x_1,\cdots,x_d)$ for the
matrix in $\m$ whose $i$th column is the column vector $x_i$. For
$r\in\{0,\cdots,d\}$ and $A,B\in\m$ we define
\bee\label{mr}
\m^r_{A,B} ~=~ \{ G\in \m: \, G_i=A_i ~\text{or}~ G_i=B_i
~\text{with}~\#\{i:\,G_i=A_i\}=r\}.
\eee
In other words,  $\m^r_{A,B}$ is the set
of all matrices $G\in\m$ with $r$ columns equal to those of $A$
(at the same places), and the remaining $d-r$ ones equal to those of $B$.
Let us define
\bee\label{1-4}
\ga_r(A,B)=\sum_{G\in\m^r_{A,B}}\det(G).
\eee

We demonstrate the main ideas for a deterministic problem first. Let
$A\in \m$ be an unknown matrix with rank $r$. Assume that, although
$A$ is unknown, we have a way of computing $\det (A+hB)$ for all $h>0$ and
some given matrix $B\in\m_d$. The multi-linearity property of the
determinant implies the following asymptotic expansion
\bee\label{detexp}
\det(A+hB) = h^{d-r} \ga_r(A,B) + O(h^{d-r+1}),
\eee
which is the core of our method.
Thus, if $\ga_r(A,B) \not = 0$, we have
\begin{align} \label{detconv}
\frac{\det (A+2hB)}{\det (A+hB)} \rightarrow 2^{d-r} \qquad \text{as}
\quad h\downarrow 0.
\end{align}
and this convergence identifies the parameter $r$. However, it is
impossible to choose a matrix $B\in \m$ which guarantees $\ga_r(A,B) \not = 0$
for all $A\in\m_r$. To solve this problem we can use a random perturbation.
As we will show later, for any $A\in\m_r$ we have $\ga_r(A,B) \not = 0$ a.s.
when $B$ is the random matrix whose entries are independent standard normal.
This idea will be the core of our testing procedure.

\subsection{Assumptions and the test statistic} \label{sec2-3}

Before we proceed with the definition of the test statistic, we introduce the
main assumptions. We need  more
structure than the mere Equation (\ref{1-0}), namely that the processes
$b_t$ and $\si_t$, and also the volatility of $\si_t$, are continuous It\^o
semimartingales. In view of the previous discussion, it is no restriction to
assume that all these are driven by the same $q$-dimensional Brownian motion,
provided we take $q$ large enough. This leads us to put
\vsq

\nib Assumption (H): \rm The $d$-dimensional semimartingale $X$, defined
on $\proba$, has the form
\bee\label{1-1}
\begin{array}{l}
X_t= X_0+ \int_0^t b_s \,ds + \int_0^t \si_s\, dW_s  \\[1.3mm]
\si_t=\si_0+ \int_0^t a_s\, ds + \int_0^t v_s\, dW_s\\[1.3mm]
b_t=b_0+ \int_0^t a'_s\, ds + \int_0^t v'_s\, dW_s \\[1.3mm]
v_t=v_0+ \int_0^t a''_s\, ds + \int_0^t v''_s\, dW_s,
\end{array}
\eee
where $W$ is a $q$-dimensional Brownian motion, and $b_t$ and $a'_t$ are
$\R^d$-valued, $\si_t$, $a_t$ and $v'_t$ are
$\R^{d\times q}$-valued, $v_t$ and $a''_t$ are $\R^{d\times q\times q}$-valued,
and $v''_t$ is $\R^{d\times q\times q\times q}$-valued, all those processes
being adapted.
Finally, the processes
$a_t,v'_t,v''_t$ are c\`adl\`ag and the processes $a'_t,a''_t$
are locally bounded.\qed
\vsq
At this stage it is not quite clear why the full force of assumption (H) is
required. In the standard limit theory for high frequency data of continuous
It\^o semimartingales, see e.g. \cite{BGJPS,J}, only
the first two representations of \eqref{1-1} are assumed. We will
further explain condition (H) once we introduce the test statistic.
When $b_t=g_1(X_t)$, $\si_t =g_2(X_t)$ with $g_1\in C^{2}(\R^d)$ and
$g_2\in C^{4}(\R^d)$, then (H) is automatically
satisfied, due to It\^o's formula.

\begin{rem}\label{R-11}  \rm
Since $\si_t$ is not uniquely specified, whereas $c_t$ is, and since
we really are interested in specific properties of $c_t$,
it would be much nicer to replace the structural assumption on $\si_t$
(second equation in (\ref{1-1})) by a similar assumption on the process
$c_t$ itself.

This is of course a trivial matter when $c_t$ is everywhere invertible:
in this case $c_t$ is a continuous It\^o semimartingale if and only if
$\si_t$ is. But here we are precisely trying to describe the rank of the
matrix $c_t$, so it is out of the question to assume that
it is {\em a priori} invertible.
Unfortunately, we were unable to replace the assumption on $\si$ by
a similar (and {\em de facto} weaker) assumption on $c$. \qed
\end{rem}

Motivated by the matrix perturbation at \eqref{detexp},
our tests will be based on statistics involving sums of (squared)
determinants. The test function will be the nonnegative map $f$
on $(\R^d)^d$ defined as
\bee\label{1-3}
f(x_1, \ldots, x_d)= \det(\text{mat}(x_1, \cdots, x_d))^2.
\eee

The authors of \cite{JLT} used the following statistics
\bee\label{I-20}
\De_n \sum_{i=1}^{[t/\De_n]-d+1}
f\big(\dd X/\rdn,\cdots,\De_{i+d-1}^n X/\rdn\big), \qquad \dd X=X_{i\De_n}-X_{(i-1)\De_n},
\eee
to test for the full rank, thus allowing for efficient
testing of the null hypothesis $\Om_T^d$. On the sets $\Om_T^r$ with $r<d$,
however, it exhibits complex degeneracies and becomes difficult to
study. In order to be able to analyze the asymptotic behavior of the preceding
statistics, we introduce a random perturbation of the original process $X$ as
motivated at the end of Subsection \ref{sec2-2} (a somewhat similar idea in a
different context was applied in \cite{CM}).
More specifically, we choose a non-random invertible $d\times d$
matrix $\Wsi$ and generate a new process
$$X'_t= \Wsi W'_t,$$
where $W'$ is a $d$-dimensional Brownian motion independent of
all processes in (\ref{1-1})(without loss of generality, for the mathematical
treatment below we may assume that it is also defined on the space $\proba$).
Following the ideas of section \ref{sec2-2}, we add to $X$ this new process
$X'$, with a multiplicative factor going to $0$. As a matter of fact we
introduce two such additions, and for $\ka=1$  or $2$ we set
\bee\label{1-6}
Z_t^{n,\ka} =X_t + \sqrt{\ka\De_n}\, X'_t.
\eee
Hence, with the notation of section \ref{sec2-2}, we use $h=\sqrt{\De_n}$,
which leads later to the optimal rate of convergence.

Another problem arises, namely in (\ref{I-20}) successive summands partly
use the same increments of $X$, and this causes problems for the Central
Limit Theorem. These problems can actually be overcome, at the expense of
quite many additional technicalities, and with the advantage of a smaller
asymptotic variance for our estimators below. However, in our case the
crucial point is the choice of the tuning ``parameter'' $\Wsi$:
this choice has an impact on the asymptotic variance
as well, and since an ``optimal'' choice of $\Wsi$ seems out of reach,
we will content ourselves with an arbitrary choice of $\Wsi$ and with a
version of (\ref{I-20}) with no overlapping of increments between the
successive summands.
This leads us to use the following two basic statistics:
\begin{align}
\label{1-7} S^{n,1}_t &= 2d\De_n \sum\limits_{i=0}^{[t/2d\De_n]-1}
f\left(\frac{Z^{n,1}_{(2id+1)\De_n}-Z^{n,1}_{2id\De_n}}{\sqrt{2\De_n}},\cdots,
\frac{Z^{n,1}_{(2id+d)\De_n}-Z^{n,1}_{(2id+d-1)\De_n}}{\sqrt{2\De_n}}\right)\\
S^{n,2}_t &= 2d\De_n \sum\limits_{i=0}^{[t/2d\De_n]-1}
f\left(\frac{Z^{n,1}_{(2id+2)\De_n}-Z^{n,1}_{(2id)\De_n}}{\sqrt{2\De_n}},\cdots,
\frac{Z^{n,1}_{(2id+2d)\De_n}-Z^{n,1}_{(2id+2d-2)\De_n}}{\sqrt{2\De_n}}\right) \nonumber.
\end{align}
Notice that the statistics $S^{n,1}_t$ and $S^{n,2}_t$ are essentially the same, except $S^{n,2}_t$
is computed using the frequency $2\De_n$.
At stage $n$ one observes the increments $\dd X$ and simulates the
increments $\dd X'$ for $i\leq[t/\De_n]$, so one ``observes'' all variables
incurring in the definition of these two statistics.

\begin{rem} \label{remstat} \rm
Now, let us explain why the assumption (H) and the random perturbation in
\eqref{1-6} are required. A direct stochastic expansion of the increments
$\De_i^n Z^{n}$ under assumption (H) implies the decomposition
\begin{align} \label{incrdec}
\text{mat}(\De_i^n Z^{n}/\sqrt{\De_n}, \ldots, \De_{i+d-1}^n Z^{n}/\sqrt{\De_n})
= \alpha_i^n + \sqrt{\De_n}(\beta_i^n(1) + \beta_i^n(2))
+ O_{\mathbb P} (\De_n),
\end{align}
where the matrices $\alpha_i^n=\text{mat}(\alpha_{i,1}^n, \ldots,
\alpha_{i,d}^n)$, $\beta_i^n(k)=
\text{mat}(\beta_{i,1}^n(k), \ldots, \beta_{i,d}^n(k))$, $k=1,2$,
in $\m$ are given by
\bee\label{incrdec-2}
\begin{array}{l}
\alpha_{i,j}^n = \De_n^{-1/2} \si_{(i-1)\Delta_n} \Delta_{i+j-1}^n W, \\[1.3 ex]
\beta_{i,j}^n(1) = b_{(i-1)\Delta_n} +
\Delta_n^{-1}\, v_{(i-1)\Delta_n} \int_{(i+j-1)\Delta_n}^{(i+j)\Delta_n}
(W_s - W_{(i+j-1)\Delta_n}) dW_s \\[1.3 ex]
\beta_{i,j}^n(2) = \De_n^{-1/2} \Wsi \Delta_{i+j-1}^n W'.
\end{array}
\eee
We remark that the matrices $\alpha_i^n, \beta_i^n(1), \beta_i^n(2)$ are
$O_{\mathbb P} (1)$. In the case $r_t\leq d-1$ for all $t$, the first order term
$\alpha_i^n$, which depends on the process $\si_t$, gives a degenerate limit
when plugged in into the statistics \eqref{1-6} or \eqref{1-7}.
Hence the second order term $\sqrt{\De_n}(\beta_i^n(1) +
\beta_i^n(2))$, which involves the processes $b_t$ and $v_t$, becomes
important. Indeed, we will see in section \ref{Res} that it affects the
limits. Furthermore, it is important to control the error of the above
decomposition, and this is done by using the last two equations in
\eqref{1-1}.

The asymptotic expansion in \eqref{incrdec} is a stochastic analogue of the
perturbation presented in \eqref{detexp} (up to an
error term) with $A=\alpha_i^n$, $B=\beta_i^n(1) + \beta_i^n(2)$ and
$h=\sqrt{\De_n}$. Under assumption (H) the term $\beta_i^n(1)$
already constitutes a random perturbation of the leading matrix $\alpha_i^n$.
However, this perturbation does not guarantee that the quantity
$\gamma_r(\alpha_i^n, \beta_i^n(1))$ defined in \eqref{1-4} does not vanish
when $\text{rank}(\si_{(i-1)\Delta_n})=r$ (which is essential for our method).
To illustrate this problem let us give a simple example. Let $d=3$, $q=1$ and
define the processes
\[
dX_t^j = \si_t^j dW_t, \qquad d\si_t^j = v_t^j dW_t, \qquad j=1,2,3.
\]
(so $W$ is a one-dimensional Brownian motion.) Then
$\text{rank}(\alpha_i^n)=1$, $\text{rank}(\beta_i^n(1))=1$, and hence
$\gamma_1(\alpha_i^n, \beta_i^n(1))=0$. The presence of the new independent
process $X'$, and thus of the term $\beta_i^n(2)$, regularizes the problem.
Indeed, we will show that $\gamma_r(\alpha_i^n, \beta_i^n(1) + \beta_i^n(2))$
does not vanish whenever $\text{rank}(\si_{(i-1)\Delta_n})=r$. Finally, the
perturbation rate $h=\sqrt{\De_n}$ in front of the process $X'$ is chosen
to achieve the best rate of convergence for the normalized versions of the
statistics $S^{n,1}_t, S^{n,2}_t$.      \qed
\end{rem}

Following the expansion \eqref{detexp} we know that that the order of $\det(\alpha_i^n + \sqrt{\De_n}(\beta_i^n(1) + \beta_i^n(2)))^2$
is increasing in $r=\text{rank}(\si_{(i-1)\Delta_n})$. Consequently, as in \eqref{detconv}, the ratio $S^{n,2}_T/S^{n,1}_T$ is expected
to identify (asymptotically) the maximal rank $R_T$. The complete asymptotic theory is presented in the next section.

\section{The asymptotic results and test for the maximal rank}
\label{Res}
\setcounter{equation}{0}
\renewcommand{\theequation}{\thesection.\arabic{equation}}

\subsection{Notation} \label{sec3-1}
In order to present the main asymptotic results we need to introduce a
few more notation. We define the function $F_r$ on $(\R^{2d})^d$ by
\bee\label{1-30}
F_r(v_1, \ldots,v_d)=
\ga_r(\text{mat}(x_1, \cdots, x_d),\text{mat}(y_1, \cdots,y_d))^2~~~\text{if}~
v_j=\left(\begin{array}{l}x_j\\y_j\end{array}\right)\in\R^{2d}.
\eee
Next, let $\ua=\m'\times\m\times\R^{dq^2}\times\R^d$, whose points are
$\uu=(\al,\be,\ga,a)$, where $\al\in\m'$ and $\be\in\m$ and
$\ga\in\R^{dq^2}$ and $a\in\R^d$. Let us denote by $\BW$ and $\BW'$ two
independent Brownian motions with respective dimensions $q$ and $d$, defined
on some space $(\BOm,\Bfa,(\Bfa_t),\BPP)$. If $\uu\in\ua$ and $\ka=1,2$ and
$i\geq1$ we associate the $2d$-dimensional variables with the following
components for $l=1,\cdots,d$:
\bee\label{1-31}
\begin{array}{l}
\Psi(\uu,\ka)^{l}_i=\frac1{\sqrt{\ka}}
\sum_{m=1}^q\al^{lm}(\BW_{\ka i}^m-\BW^m_{\ka(i-1)})\\[1.5 ex]
\Psi(\uu,\ka)^{d+l}_i=a^l+\frac1{\sqrt{\ka}}
\sum_{m=1}^d\be^{lm}(\BW_{\ka i}'^m-\BW'^m_{\ka(i-1)})
+\frac1{\ka}\sum_{m,k=1}^q\ga^{lmk}
\int_{\ka(i-1)}^{\ka i}\BW_s^k\,d\BW_s^m.\end{array}
\eee
With the notation \eqref{1-30} we can then define the variables
\bee\label{1-32}
\BF_r(\uu,\ka)=
F_r\big(\Psi(\uu,\ka)_1,\cdots,\Psi(\uu,\ka)_d\big).
\eee
The two sequences $(\Psi(\uu,\ka))_{i\geq1}$ are not independent, but they
have the same (global) law, for $\ka=1,2$. Therefore if $\uu=(\al,\be,\ga,a)$
we can set
\bee\label{1-8}
\begin{array}{l}
\Ga_r(\uu)=\BE\big(\BF_r(\uu,1)\big)=\BE\big(\BF_r(\uu,2)\big)\\[1.5 ex]
\Ga'_r(\uu)=\BE\big(\BF_r(\uu,1)^2\big)-\Ga_r(\uu)^2
=\BE\big(\BF_r(\uu,2)^2\big)-\Ga_r(\uu)^2\\[1.5 ex]
\Ga''_r(\uu)=\BE\big(\BF_r(\uu,1)\,\BF_r(\uu,2)\big)-\Ga_r(\uu)^2.
\end{array}
\eee
We then obtain the following crucial properties

\begin{lem}\label{L-0}
Let $\uu=(\al,\be,\ga,a)\in\ua$ with $\be\in\m_d$. Then if $r\in
\{0,1,\cdots,d\}$,
\begin{align}
\label{1-9} \text{\rm rank}(\al)=r~&\Longrightarrow~\Ga_r(\uu)>0,\quad
\Ga_r'(\uu)>\Ga_r''(\uu)   \\[1.5 ex]
\text{\rm rank}(\al)<r~&\Longrightarrow~
\Ga_r(\uu)=\Ga'_r(\uu)=\Ga''_r(\uu)=0. \nonumber
\end{align}
\end{lem}

\subsection{The limiting results} \label{sec3-2}

The key result is the asymptotic behavior of the processes
$S^{n,j}$ as $n\to\infty$. These processes enjoy a Law of Large Numbers
and a Central Limit Theorem, the centering being around one of the following
processes, where $r$ is any (fixed) integer between $0$ and $r$:
\bee\label{2-2}
S(r)_t~=~\int_0^t\Ga_r(\si_s,\Wsi,v_s,b_s)\,ds\,.
\eee
We will in fact have a CLT for the two-dimensional processes
$U(r)^n$ with components
\bee\label{2-101}
U(r)^{n,\ka}=\rdnn\,\Big(\frac1{(\ka\De_n)^{d-r}}\,S^{n,\ka}-S(r)\Big).
\eee
Of course, the centering process $S(r)$ depends on $r$, so one needs an
additional assumption related with the particular value of $r$ which is
chosen below (in contrast, the centering term is the same for all
components):

\begin{theo}\label{T1} Assume (H), and also that $r_t(\om)\leq r$ identically
for some $r\in\{0,\cdots,d\}$.  Then we have the stable (functional)
convergence in law
\bee\label{2-1}
U(r)^n~\tolls~\ua(r),
\eee
where $\ua(r)=(\ua(r)^\ka)_{\ka=1,2}$ is defined on an extension $\probt$ of
$\proba$ and is, conditionally on $\f$, a continuous centered Gaussian
martingale with conditional covariance
\bee\label{2-5}
\WE(\ua(r)^\ka_t\,\ua(r)^{\ka'}_t\mid\f)=V(r)^{\ka\ka'}_t:=
\left\{\begin{array}{ll}
2d\int_0^t\Ga_r'(\si_s,\Wsi,v_s,b_s)\,ds&\text{\rm if}~\ka=\ka'\\[1.5 ex]
2d\int_0^t\Ga_r''(\si_s,\Wsi,v_s,b_s)\,ds&\text{\rm if}~\ka\neq \ka'.
\end{array}\right.
\eee
\end{theo}
Note that in the above setting,  if $r<r'\leq d$, we also have $r_t\leq r'$
and thus the results also hold with $r'$ instead of $r$ everywhere. This does
not bring a contradiction because, by (\ref{1-9}), in this case the processes
$S(r')$ and $U(r')$ are identically vanishing.

Now, these processes $S^{n,j}$ are only tools, and at the end we
will be interested, for any $T>0$ fixed, in ``estimators'' for $R_T$, which
are
\bee\label{2-301}
\wR(n,T)=d-\frac{\log(S^{n,2}_T/S^{n,1}_T)}{\log2}.
\eee
The quantity is a transformed analogue of the term on the left side of
\eqref{detconv}.
The following corollary is then a simple consequence of the previous theorem:

\begin{cor}\label{C12} Assume (H), and let $r\in\{0,\cdots,d\}$ and $T>0$.
Then the following stable convergence in law holds:
\bee\label{2-6}
\rdnn\,(\wR(n,T)-r)~\tols~\s(T)\quad
\quad \text{\rm on the set}~\Om^r_T,
\eee
where $\s(T)$ can be realized as $\s(T)=\frac1{\log2}\,
(\ua(r)^1_T-\ua(r)^2_T)/S(r)_T$
and is thus defined on an extension $\probt$ of $\proba$
and is, conditionally on $\f$, a centered Gaussian variable
whose conditional variance is
$$\WE((\s(T))^2\mid\f)=V(T),$$
where $V(T)$ is a.s. positive and given by
\bee\label{2-7}
V(T)=\frac1{(\log2)^2}\,\frac{V(r)^{1,1}_T+V(r)_T^{2,2}
-2V(r)^{1,2}_T}{(S(r)_T)^2}
\quad\text{\rm on each set $\Om^r_T$}.
\eee
\end{cor}
In order to make this result feasible, we need consistent estimators for
$V(T)$. For the denominator $S(r)_T^2$ we
can of course take the square of $\De_n^{r-d}\,S_T^{n,1}$. As for the
numerator, we need estimators for $V(r)^{\ka,\ka'}_T$. Up
to normalization, natural ones are as follows:
\begin{align}
 V_t^{n,\ka\ka'}=4d^2\De_n \sum\limits_{i=0}^{[t/2d\De_n]-1}
f\left(\frac{Z^{n,\ka}_{(2id+\ka)\De_n}-Z^{n,\ka}_{2id\De_n}}{\sqrt{\ka\De_n}},
\cdots,
\frac{Z^{n,\ka}_{(2id+\ka d)\De_n}-Z^{n,\ka}_{(2id+\ka(d-1))\De_n}}
{\sqrt{\ka\De_n}}\right) \nonumber \\\hskip4cm
\label{2-26} \times
f\left(\frac{Z^{n,\ka'}_{(2id+\ka')\De_n}-Z^{n,\ka'}_{2id\De_n}}
{\sqrt{\ka'\De_n}},\cdots,
\frac{Z^{n,\ka'}_{(2id+\ka' d)\De_n}-Z^{n,\ka'}_{(2id+\ka'(d-1))\De_n}}
{\sqrt{\ka'\De_n}}\right).
\end{align}

\begin{prop}\label{P3} Assume (H).

a) If $r_t(\om)\leq r$ identically
for some $r\in\{0,\cdots,d\}$, we have for $\ka,\ka'=1,2$:
\bee\label{2-29}
\begin{array}{c}
\frac1{(\ka\ka'\De_n^2)^{d-r}}\,V^{n,\ka\ka'}~\toucp~
2d\int_0^\cdot\Te^{r,\ka,\ka'}_s\,ds,\qquad\text{\rm where}\\[1.5 ex]
\Te^{r,\ka,\ka'}_s=\left\{\begin{array}{ll}
\Ga'_r(\si_s,\Wsi,v_s,b_s)+\Ga_r(\si_s,\Wsi,v_s,b_s)^2
~&\text{\rm if}~\ka=\ka'\\[1.5mm]
\Ga''_r(\si_s,\Wsi,v_s,b_s)+\Ga_r(\si_s,\Wsi,v_s,b_s)^2
~&\text{\rm if}~\ka\neq\ka'.\end{array}\right.
\end{array}
\eee

b) We have
\bee\label{2-30}
V(n,T):=\frac{V^{n,11}_T+2^{2(\wR(n,T)-d)}V^{n,22}_T
-2^{1+\wR(n,T)-d}V^{n,12}_T}{(S^{n,1}_T\,\log2)^2}~\toop~ V(T).
\eee
\end{prop}

\begin{rem}\label{R1} \rm The numerator of the right side of
\eqref{2-7} is also $2(V(r)^{11}_T-V(r)^{12}_T)$. Therefore we have
$$V'(n,T)=\frac{V^{n,11}_T
-2^{1+\wR(n,T)-d}V^{n,12}_T}{(S^{n,1}_T)^2}~\toop~ V(T)\quad
\text{on the set}~\Om_T^r$$
as well. However, $V(n,T)\geq 0$ by construction (and it is even a.s.
positive unless $r_t=0$ identically on $[0,T]$), a property not
shared by $V'(n,T)$.\qed
\end{rem}

Now, by the delta-method for stable convergence
in law, the two previous results immediately yield:

\begin{cor}\label{C3} Under (H) and for any $T>0$ we have
\bee\label{2-67}
\frac{\wR(n,T)-R_T}{\sqrt{\De_n\,V(n,T)}}~\tols~\Phi,
\eee
where $\Phi\sim \mathcal N(0,1)$ is defined on an extension $\probt$ of $\proba$ and is independent of $\mathcal F$.
\end{cor}

\subsection{Tests for the maximal rank} \label{sec3-3}
So far, it seems that $\wR(n,T)$ are estimators for the maximal rank
$R_T$, which equals $r$ on the set $\Om^r_T$, and even feasible estimators
if we use Corollary \ref{C3}. In particular, this corollary seems to
allow us to easily construct confidence intervals for $R_T$.

However, what precedes does not make much statistical sense: the parameter
$R_T$ to be estimated takes its values in $\{0,1,\cdots,d\}$, whereas the
estimators $\wR(n,T)$ are of course not integer-valued and can even be
negative, or bigger than $d$. One could overcome this problem by taking the
integer closest to $\wR(n,T)$, say
$\wR'(n,T)$, and then use $\wR''(n,Y)=0\vee(\wR'(n,T)\wedge d)$ as the final
estimator. Note that $\wR''(n,T)$ enjoys the same CLT as $\wR(n,T)$ does, on
each $\Om^r_T$ with $1\leq r\leq d-1$, but of course not when $R_T=0$ or
$R_T=d$, in which cases the limiting law of the normalized error is ``half
Gaussian and half a Dirac mass at $0$''.
Furthermore, confidence intervals have little meaning in this context, except
perhaps when the dimension of $X$ is very large.

So, it seems more appropriate here to do testing: we can test the null
hypothesis
that the path lies in $\Om_T^r$ for some $r$, against the alternative
that it is in $\Om_T^{r'}$ for another specific $r'\neq r$, or for
all $r'>r$ or all $r'<r$, or all $r'\neq r$. We may also use composite
null hypotheses, such as being in $\Om_T^r$ for some $r$ smaller, or bigger,
than a given value $r_0$.

We start with the problem of testing the null hypothesis $\Om_T^r$,
against the alternative $\Om_T^{\neq r}=\cup_{r'\neq r,0\leq r'\leq d}
\Om_T^{r'}$. For any $\al\in(0,1)$, and with $z_\al$ being the symmetric
$\al$-quantile of $\n(0,1)$ defined by $\PP(|\Phi|>z_\al)=\al$ when
$\Phi\sim\n(0,1)$, we take the critical (rejection) region
\bee\label{2-60}
\ca(\al)_T^{n,=r}=\big\{\om:~|\wR(n,T)-r|>z_\al\sqrt{\De_n\,V(n,T)}
\big\}.
\eee

\begin{prop}\label{P1} Under (H), the tests (\ref{2-60}) have the
asymptotic level $\al$ for testing the null $\Om_T^r$, in the sense that
\bee\label{2-61}
A\subset\Om_T^r,~\PP(A)>0~~\Rightarrow~~
\PP\big(\ca(\al)_T^{n,=r}\mid A\big)\to\al
\eee
(above, $\PP(.\mid A)$ is the usual conditional probability). They are also
consistent for the alternative $\Om_T^{\neq r}$, in the sense that
\bee\label{2-62}
\PP\big(\ca(\al)_T^{n,=r}\cap\Om_T^{\neq r}\big)\to\PP(\Om_T^{\neq r}).
\eee
\end{prop}

One constructs one-sided tests in the same way. For example, if we want
to test the null hypothesis $\Om_T^{\leq r}=\cup_{r'\leq r}\Om_T^{r'}$
against the alternative $\Om_T^{>r}=\cup_{r'>r}\Om_T^{r'}$, and if $z'_\al$
is the one-sided $\al$-quantile defined
by $\PP(\Phi>z'_\al)=\al$, we take the critical region

\bee\label{2-63}
\ca(\al)_T^{n,\geq r}=\big\{\om:~\wR(n,T)>r+z'_\al\sqrt{\De_n\,V(n,T)}
\big\}.
\eee
Exactly as above, one obtains the following proposition.

\begin{prop}\label{P2} Under (H), the tests (\ref{2-63}) have the
asymptotic level at most $\al$ for testing the null $\Om_T^{\leq r}$, and
indeed satisfy
\bee\label{2-64}
A\subset\Om_T^{\geq r},~\PP(A)>0~~\Rightarrow~~
\PP\big(\ca(\al)_T^{n,=r}\mid A\big)\to\al\,\PP(\Om_T^r\mid A)\leq\al,
\eee
and are consistent for the alternative $\Om_T^{>r}$.
\end{prop}

The tests for the null $\Om_T^{\geq r}$ against $\Om_T^{<r}$ are obtained
analogously.

\begin{rem} \label{relprob} \rm
Let us link our testing procedure with some other statistical problems:
\vsd

\ni a) In \cite{DS, G} parametric estimation methods for the so called integrated
diffusions have been developed. An integrated diffusion is a
process that satisfies the first and the third equations of assumption (H)
with $\si=0$, i.e.
\[
dX_t= b_t dt,
\]
where $b_t$ is a continuous It\^o semimartingale. We refer to \cite{DS} for
various applications of these models in natural sciences.
Given high frequency observations of $X$, testing the null hypothesis of
integrated diffusion versus the alternative of a diffusion
with a present volatility part $\sigma$ is equivalent to testing $\Om_T^0$
versus $\Om_T^{\neq 0}$.
\vsd

\ni b) Another potential application of our
method is a test for ``perfect correlation'' between the process $X$ and the
unobserved volatility $\si$. The problem can be formulated as follows: Let $X$
and $\sigma$ be two one-dimensional continuous It\^o semimartingales
of the form
\[
dX_t= b_t dt + \sigma_t dW_t, \qquad d\sigma_t^2= a_t dt + v_t dB_t,
\]
where $W$ and $B$ are one-dimensional Brownian motions with the bracket
process $[W,B]_t=\rho t$, $|\rho|\leq 1$. For financial
applications testing the hypothesis $|\rho|= 1$ versus $|\rho|< 1$ is of
certain interest. Note that $|\rho|= 1$ appears
in the SDE case, i.e. when $\si_t=g(X_t)$ with $g\in C^2(\R)$.  We refer to
testing local volatility hypothesis
in \cite{PR} for a more detailed discussion (see also \cite{V} for related
statistical problems). The aforementioned problem
is equivalent to testing $\Om_T^1$ versus $\Om_T^{> 1}=\Om_T^2$ for the
two-dimensional process $(X,\si^2)$. Since the process
$\si^2$ is unobserved, it has to be locally estimated from the high frequency
observations of $X$ first (see e.g. \cite{PR} for more details). \qed
\end{rem}

\section{A test for a constant rank} \label{sec4}

This section is devoted to a seemingly different topic, namely
whether the {\em a priori} time-dependent rank is constant or not.
Our test statistics will be based on a distance measure between the rank
process $r_t$ and the maximal rank $R_T$, which vanishes if and only if the
rank is constant almost surely. For the formal testing procedure we will need
some limiting results for the ``spot estimators''
of the rank. By this, we mean estimators for $r_t$, for any given $t$,
at least under the assumption that $r_s$ is equal to $r_t$ for all $s$
in some right or left neighborhood of $t$.

To describe these spot estimators we pick a sequence $k_n\geq1$ of integers
going to infinity, and such that $k_n\De_n\to0$ (as for spot volatility
estimators), and precise specifications for $k_n$ will be given later,
although we always assume $k_n\geq4d$. For any integer $i\geq1$ we
set
\bee\label{2-68}
\wR^n_i=d-\frac{\log \wS^n_i}{\log2},\qquad
\wS^n_i=\frac{S^{n,2}_{2d(i+1)k_n\De_n}-S^{n,2}_{2dik_n\De_n}}
{S^{n,1}_{2d(i+1)k_n\De_n}-S^{n,1}_{2dik_n\De_n}}.
\eee
Then $\wR^n_i$, more or less, plays the role of an estimator of the maximum
of $r_t$ over an interval of length $2dk_n\De_n$ around the time
$2id\De_n$, and we set for any $p>0$:
\bee\label{2-69}
\begin{array}{c}
A(p)^n_t=2dk_n\De_n\sum_{i=0}^{[t/2dk_n\De_n]-2}
\{|\wR^n_{ik_n}|^p\wedge(d+1)^p\}\\[1.5 ex]
\!B(n,p,T)=A(p)^n_T-a(n,T)(\wR(n,T))^p,\quad
a(n,T)=2dk_n\De_n\,\big([T/2dk_n\De_n]-1\big).\end{array}
\eee
The asymptotic results for the quantity $B(n,p,T)$ are as follows.

\begin{theo}\label{T4} Assume (H), and let $T>0$, $p>0$ and $k_n$ be such that
$k_n\De_n^{3/4}\to\infty$ and $k_n\De_n\to0$.

a) If $t\mapsto r_t(\om)$ is continuous except at finitely many points on
$[0,T]$, hence piecewise constant, we have
\bee\label{2-70}
B(n,p,T)~\toop~\int_0^T(r_s)^p\,ds-T(R_T)^p.
\eee

b) We have the stable convergence in law:
\bee\label{2-71}
\rdnn\,B(n,p,T)~\tols~ \ba(p,T)\quad\text{\rm in restriction to the set}~
\Om_T^=\cap\{R_T\geq1\},
\eee
where $\ba(p,T)$ is defined on an extension $\probt$ of $\proba$ and is,
conditionally on $\f$, a centered Gaussian variable with conditional variance
$\BV(p,T)=\WE(\ba(p,T)^2\mid\f)$ given on each set $\Om_T^r$ by
\bee\label{2-72}
\BV(p,T)=\left(\frac{pr^{p-1}}{\log2}\right)^2
\int_0^T\left(\frac1{\Ga_r(\si_s,\Wsi,v_s,b_s)}-\frac T{S(r)_T}\right)^2
\,\big(dV(r)^{11}_s+dV(r)^{22}_s-2dV(r)^{12}_s\big),
\eee
with $V(r)^{\ka \ka'}$ being defined at \eqref{2-5}.
\end{theo}

Notice that the right side of \eqref{2-70} is $0$ on the set $\Om_T^=$, and
strictly negative on $\Om_T^{\neq}$.

\begin{rem}\label{R-14} \rm The reader will notice that in the definition
of $A(p)^n_t$ the summands are $|\wR^n_i|^p\wedge(d+1)^d$, instead of
the more natural $|\wR^n_i|^p$. We could take this more natural form for
(b) above, but it is useful (and innocuous from a practical viewpoint)
to ``bound'' the summands, in order to obtain (a). We could bound them
by $d^p$ instead of $(d+1)^p$ and still have \eqref{2-70}, but then
\eqref{2-71} would then fail in case $r=d$ is the maximal rank: we would
obtain a CLT with a non-Gaussian and non-centered limit.\qed
\end{rem}

\begin{rem}\label{R-15} \rm In the setting of (b) above, we will in fact
prove a joint convergence for the variables $A(p)^n_T-a(n,T)r^p$ and
$\wR(n,T)-r$, both normalized by $1/\rdn$ (the second one being as in
\eqref{2-6}), and from which \eqref{2-71} follows. Such a joint CLT even
holds under the assumptions of (a), with a complicated limit, but this
refinement is not useful for us in this paper.
\end{rem}

\begin{rem}\label{R-16} \rm
One can also prove a joint convergence for the variables $A(p)^n_T-a(n,T)r^p$
with different values of $p$, and still normalized by $1/\rdn$. However,
when $p>p'>0$ it turn out that the difference
$\rdnn\,(A(p)^n_T-a(n,T)^{1-p'/p}\,(A(p')^n_T)^{p/p'})$ converges to $0$, and
no known normalization gives a proper CLT.\qed
\end{rem}

As before, we need consistent estimators for the conditional variance
$\BV(p,T)$. Such estimators are constructed in a way analogous to
\eqref{2-26}. That is, we set with $k_n$ as above:
\begin{align}
 \BV_t^{n,\ka\ka'}&=4d^2\De_n^{1+2d-2\wR(n,T)}
\sum\limits_{i=0}^{[t/2d\De_n]-k_n-1}
\Big(\frac{2dk_n\De_n}{S^{n,1}_{2d(i+k_n)\De_n}
-S^{n,1}_{2id\De_n}}-\frac{T}{S^{n,1}_T}\Big)^2 \nonumber\\
&\times
f\Big(\frac{Z^{n,\ka}_{(2id+\ka)\De_n}-Z^{n,\ka}_{2id\De_n}}{\sqrt{\ka\De_n}},
\cdots,
\frac{Z^{n,\ka}_{(2id+\ka d)\De_n}-Z^{n,\ka}_{(2id+\ka(d-1))\De_n}}
{\sqrt{\ka\De_n}}\Big) \label{2-50} \\
&\times
f\Big(\frac{Z^{n,\ka'}_{(2id+\ka')\De_n}-Z^{n,\ka'}_{2id\De_n}}
{\sqrt{\ka'\De_n}},\cdots,
\frac{Z^{n,\ka'}_{(2id+\ka' d)\De_n}-Z^{n,\ka'}_{(2id+\ka'(d-1))\De_n}}
{\sqrt{\ka'\De_n}}\Big). \nonumber
\end{align}

\begin{theo}\label{T5} Assume (H), and let $T>0$, $p>0$ and $k_n$ be such that
$k_n\De_n^{3/4}\to\infty$ and $k_n\De_n\to0$. Then we have
\bee\label{2-52}
\begin{array}{l}
\BV(n,p,T):=\Big(\frac{p\wR(n,T)^{p-1}}{\log2}\Big)^2
\,\big(\BV^{n,11}_T+2^{2(\wR(n,T)-d)}\BV^{n,22}_T
-2^{1+\wR(n,T)-d}\BV^{n,12}_T\big)\\
\hskip5cm ~\toop~ \BV(p,T)\quad\text{\rm on the set}~\Om^=_T.
\end{array}
\eee
Moreover, the variables
\bee\label{2-53}
Z(n,p,T)=\frac{B(n,p,T)}{\sqrt{\De_n\,(\BV(n,p,T)\wedge(1/\rdn\,))}}
\eee
have the following asymptotic behavior, where $\Phi\sim \mathcal N(0,1)$ is as in Corollary \ref{C3}:
\bee\label{2-54}
\begin{array}{ll}
Z(n,p,T)~\tols~\Phi~&\text{\rm in restriction to the set}~
\Om_T^=\cap\{R_T\geq1\}\\
Z(n,p,T)~\toop~-\infty~&\text{\rm in restriction to the set}~\Om_T^{\neq}
\end{array}
\eee
\end{theo}

Having all instruments at hand we proceed with testing. What
is easily available is a family of tests for the null $\Om_T^=$, whereas
the alternative is restricted to $\Om_T^{\neq}$. One does not know how
to test the null $\Om_T^{\neq}$.

For this purpose we use the statistic $B(n,p,T)$. In fact, \eqref{2-54} gives us
the behavior of this statistic on $\Om_T^=\cap\{R_T\geq1\}$, and this
is the null which is tested below. Now, $\Om^=_T$ is the union
of $\Om_T^=\cap\{R_T\geq1\}$ and $\Om^0_T$, so if we are interested in
testing the whole $\Om^=_T$ one can do a double test, using what
precedes and Proposition \ref{P2} with $r=0$.

We propose to use the following critical region, where $p>0$ is chosen
arbitrarily and $z'_\al$ is again the one-sided $\al$-quantile of
$\n(0,1)$:
\bee\label{2-65}
\ca(\al)_T^{n,\equiv}=\big\{\om:~B(n,p,T)<-
z'_\al\sqrt{\De_n\,(\BV(n,p,T)\wedge(1/\rdn\,))}\big\}.
\eee
Exactly as in the previous section we obtain the following result.

\begin{prop}\label{P20} Under (H), the tests (\ref{2-65}) have the
asymptotic level $\al$ for testing the null $\Om_T^=\cap\{R_T\geq1\}$,
in the sense of (\ref{2-61}), and
are consistent for the alternative $\Om_T^{\neq}$.
\end{prop}

\section{Proofs}\label{P}
\setcounter{equation}{0}
\renewcommand{\theequation}{\thesection.\arabic{equation}}

Before we start presenting the formal proofs, let us give the road map.
Subsection \ref{sec5-1} demonstrates some technical results on expansions of
determinants. They are applied in Subsection \ref{sec5-2} to prove Lemma
\ref{L-0}. This Lemma implies that the process $S(r)_t$ defined at
\eqref{2-2} is strictly positive on the set $\Om^r_T$, which is crucial
for our method.

The first main result of our paper is Theorem \ref{T1} whose proof is rather
involved. First, we will show that the standard localization
procedure (see e.g. Section 3 in \cite{BGJPS}) implies that all processes in
(H) may be assumed to be bounded without loss of generality. This first
step considerably simplifies the stochastic treatment of various quantities.
A second crucial step is the stochastic expansion explained in
Remark \ref{remstat}: we have \eqref{incrdec} and \eqref{incrdec-2}.
Subsection \ref{sec5-3} deals with the formal justification of this expansion,
for which we will use slightly different notation.

It turns out that the stochastic order of the error term related to the
decomposition \eqref{incrdec}, namely $O_{\mathbb P} (\De_n)$, is not
sufficient to show its asymptotic negligibility. However, we will prove that
the error terms are martingale differences, so they will not affect the stable
central limit theorem at \eqref{2-1}. A similar treatment will be required
for the error term
connected with the stochastic version of the expansion \eqref{detexp}.

The proof of Proposition \ref{P3} (consistent estimation of the asymptotic
conditional covariance matrix) is somewhat easier. Corollary
\ref{C12} follows essentially from Theorem \ref{T1} by the delta method for
stable convergence. The proofs of these results are collected in Subsection
\ref{sec5-4}. In particular, we apply a stable central limit theorem for
semimartingales (see e.g. \cite[Theorem IX.7.28]{JS}) to prove Theorem \ref{T1}.

The proof of Theorems \ref{T4} and \ref{T5}, which is presented in Subsection
\ref{sec5-5}, is a bit more involved than one of Theorem \ref{T1}, although
the main techniques are similar. The additional difficulty comes from the fact
that we need to use the stable convergence of Theorem \ref{T1},
but for processes evaluated at random times. Corollary \ref{C3} and
Propositions \ref{P1}, \ref{P2} and \ref{P20}
are straightforward consequences of the previous results.

\subsection{Expansion of determinants.} \label{sec5-1}
We first prove some general and easy facts about determinants. Below
$\|A\|$ denotes the Euclidean norm of a matrix $A\in\m$.

For $m\geq1$ we call $\p_m$ the set of all multi-integers
$\ppp=(p_1,\cdots,p_m)$ with $p_1+\cdots+p_m=d$, and
$\ia_{\ppp}$ is the set of all partitions
$\III=(I_1,\cdots,I_m)$ of $\{1,\cdots,d\}$ such that $I_j$ contains
exactly $p_j$ points (so $I_j=\emptyset$ if $p_j=0$). If $\ppp\in\p_m$ and
$\III\in\ia_\ppp$ and $A_1,\cdots,A_m\in\m$, we write
$G^\III_{A_1,\cdots,A_m}$ for the matrix
whose $i$th column is the $i$th column of $A_j$ when $i\in I_j$.
Letting $A,B,C\in\m$, we can rewrite (\ref{1-4}) as
\bee\label{P-100}
\ga_r(A,B)=~\sum_{\III\in\ia_{(r,d-r)}}\det(G^\III_{A,B}),
\eee
and we set
\bee\label{P-104}
\ga'_r(A,B,C)=\sum_{\III\in\ia_{(r,d-r-1,1)}}\det(G^\III_{A,B,C}).
\eee
In the following two lemmas we present some technical results on determinant
expansions.

\begin{lem}\label{LP-11} For any $m\geq1$ and $A_1,\cdots,A_m\in\m$ we have
\bee\label{P-101}
\det(A_1+\cdots+A_m)=\sum_{\ppp\in\p_m}~\sum_{\III\in\ia_\ppp}
\det(G^\III_{A_1,\cdots,A_m}).
\eee
\end{lem}

\nib Proof. \rm Letting $\s_d$ be the set of all permutations of
$\{1,\cdots,d\}$ and $\sign(s)$ be the signature of $s\in\s$, we have
\bean
\det(A+B)&=&\sum_{s\in\s_d}(-1)^{\sign(s)}\prod_{i=1}^d
(a^{s(i),i}+b^{s(i),i}) \\
&=&\sum_{I\subset\{1,\cdots,d\}}~
\sum_{s\in\s_d}(-1)^{\sign(s)}\prod_{i\in I}a^{s(i),i}
\prod_{i\notin I}b^{s(i),i}=\sum_{I\subset\{1,\cdots,d\}}
\det(G_{A,B}^{(I,I^c)}).
\eean
This readily implies that if (\ref{P-101}) holds for some $m$, it also holds
for $m+1$. Since (\ref{P-101}) is obvious for $m=1$, the result follows
by induction on $m$.\qed

\begin{lem}\label{LP-2} There is a constant $K$ such that, for all
$r=0,\cdots,d$, all $h\in(0,1]$ and all $A,B,C,D\in\m$ with
{\rm rank}$(A)\leq r$
we have, with $\La=\|A\|+\|B\|+\|C\|+\|D\|$ and with the convention
$\ga_{-1}(A,B)=0$:
\bee\label{P-103}
\begin{array}{l}
\!\!\big|\det(A+hB+h^2C+h^2D)-h^{d-r}\ga_r(A,B)-
h^{d-r+1}(\ga_{r-1}(A,B)+\ga'_r(A,B,C))\big|\\
\hskip4cm \leq Kh^{r-d+1}\La^{d-1}(h\La+\|D\|),
\end{array}
\eee
\bee\label{P-105}
\begin{array}{l}
\big|\frac1{h^{2d-2r}}\,\det(A+hB+h^2C+h^2D)^2-\ga_r(A,B)^2\\
\hskip2cm
-2h\,\ga_r(A,B)(\ga_{r-1}(A,B)+\ga'_r(A,B,C))\big|
\leq Kh\La^{2d-1}(h\La+\|D\|).
\end{array}
\eee
\end{lem}

\nib Proof. \rm Let $\ppp\in\p_4$ and $\III\in\ia_\ppp$. Then
$\det(G^\III_{A,hB,h^2C,h^2D})=h^{p_2+2p_3+2p_3}\det(G^\III_{A,B,C,D)})$
vanishes when $p_1>r$, and has absolute value smaller
than $Kh^{p_2+2p_3+2p_4}\La^{d-p_4}\|D\|^{p_4}$. Then (\ref{P-103})
readily follows from (\ref{P-101}), and by taking squares in
(\ref{P-103}) we deduce (\ref{P-105}).\qed \newline \newline
With the same notation, and if further $A',B',C',D'\in\m$ with
rank$(A')\leq r$ also and $\La'=\|A'\|+\|B'\|+\|C'\|+\|D'\|$, and
$h'\in(0,1]$, the same argument shows that
\bee\label{P-106}
\begin{array}{l}
\big|\frac1{(hh')^{2d-2r}}\,\det(A+hB+h^2C+h^2D)^2\,
\det(A'+h'B'+h'^2C'+h'^2D')^2\\[1.5 ex]
\hskip4cm-\ga_r(A,B)^2\,\ga_r(A',B')^2\big|\leq K(h+h')(\La\,\La')^{2d}.
\end{array}
\eee

\subsection{Proof of Lemma \ref{L-0}.} \label{sec5-2}

\nib1) The results about $\Ga_r(\uu)$. \rm
We write $V_i$ and $\BV_i$ for the $d$-dimensional
variables whose components are respectively the $d$ first and the $d$ last
components of $\Psi(\uu,1)_i$, for which we can take $\BW=W$ and $\BW'=W'$,
and we set $A=$~mat$(V_1,\cdots,V_d)$ and
$B=$~mat$(\BV_1,\cdots,\BV_d)$. If $\De_jW^{(\prime)}=
W^{(\prime)}_j-W^{(\prime)}_{j-1}$, we have
\bee\label{P-20}
V_i^l=\sum_{m=1}^q\al^{lm}\De_iW^m,\quad
\BV_i^l=a^l+\sum_{m=1}^d\be^{lm}\De_iW'^m
+\sum_{m,k=1}^q\ga^{lkm}h_{i,km}(W),
\eee
where each $h_{i,lm}$ is a function of the path of $W$. Note also
that $\BF_r(\uu,1)=\ga_r(A,B)^2$.

Assuming first that the rank
of $\al$ is (strictly) smaller than $r$, we observe that the rank of $A$
is also smaller than $r$, implying by \eqref{P-100} that
$\ga_r(A,B)=0$, hence $\Ga_r(\uu)=0$.

Next we assume that the rank of $\al$ is $r$, and proceed to prove
$\Ga_r(\uu)>0$.
We first simplify the problem as follows. The matrix $\be$ is invertible
and the rank of $\be^{-1}\al\al^*\be^{-1,*}$ is $r$, so we can write
$\be^{-1}\al=\Pi\La$,
where $\Pi\in\m$ is an orthonormal matrix and $\La\in\m$ is a diagonal matrix
whose diagonal entries $\la_j$ satisfy $\la_j\neq0$ if $j\leq r$ and
$\la_j=0$ otherwise. Then, setting $V'_j=\Pi^*\be^{-1}V_j$ and
$\BV'_j=\Pi^*\be^{-1}\BV_j$, the sequence $(V'_j,\BV'_j)$ has the form
(\ref{P-20}), upon replacing $W'$ by $\Pi^*W'$ (another Brownian motion)
and $\uu=(\al,\be,\ga,a)$ by $\uu'=(\al',J,\ga',a')$, where $J$ is the
identity in $\m$ and $\al'=\Pi^*\be^{-1}\al=\La$
and $\ga'^{ijl}=\sum_{m=1}^d(\Pi^*\be^{-1})^{jm}\ga^{mjl}$ and
$a'=\Pi^*\be^{-1}a$. Furthermore, $A'=$ mat$(V'_1,\cdots,V'_d)=\Pi^*\be^{-1}A$
and $B'=$ mat$(\BV'_1,\cdots,\BV'_d)=\Pi^*\be^{-1}B$, implying
$\ga_r(A',B')=\det(\Pi^*\be^{-1})\ga_r(A,B)$, which in turn yields
$\Ga_r(\uu)=\frac1{\det(\be)^2}\,\Ga_r(\uu')$, because $\det(\Pi)=1$.

In other words, it is enough to prove the result when $\be=J$ and
$\al=\La$ is diagonal as above, and below we assume this. The two matrices
$A$ and $B$ can thus be realized as
\bee\label{P-21}
A=\left(\begin{array}{ccc}\la_1\Phi_1^1&\cdots&
\la_1\Phi_d^1\\\cdot&\cdots&\cdot\\\cdot&\cdots&\cdot\\
\la_d\Phi_1^d&\cdots&
\la_d\Phi_d^d\end{array}\right),\qquad
B=\left(\begin{array}{ccc}\Up_1^1+\Te^1_1&\cdots&
\Up_d^1+\Te_d^1\\\cdot&\cdots&\cdot\\\cdot&\cdots&\cdot\\
\Up_1^d+\Te^d_1&\cdots&
\Up_d^d+\Te_d^d\end{array}\right),
\eee
where all $\Phi_j^i$ and $\Up_j^i$ are i.i.d. $\n(0,1)$ and the variables
$\Te^i_j$ are independent of the $\Up^l_m$'s (note that we have
incorporated the constant $a^i$ in each variable $\Te^i_j$).

Let $\ja_r$ be the class of all subsets of $\{1,\cdots,d\}$ with $r$
points. Since $\la_j\neq0$ if $j\leq r$ and $\la_j=0$ otherwise, we see
that, if $\III=(I,I^c)$ with $I=\{j_1<\cdots<j_r\}\in\ja_r$ and
$I^c=\{j'_1<\cdots<j'_{d-r}\}$, we
have $\det(G^\III_{A,B})=\ep_I\det(A_I)\,\det(B_I)$, where $A_I$ and
$B_I$ are the $r\times r$ and $(d-r)\times(d-r)$ matrices with entries
$A_I^{l,m}=A^{j_l,m}$ and $B_I^{l,m}=B^{j'_l,r+l}$, and $\ep_I$ takes values in $\{-1,1\}$. Thus
$$\ga_r(A,B)=\sum_{I\in\ja_r}\ep_I \det(A_I)\,\det(B_I).$$
In this sum we single out the $I$'s which contain $d$, and those which
do not, and for the former ones the product $\det(A_I)\det(B_I)$ does not
depend of the vector $\Up_d$. For those which do not contain
$d$, we develop $\det(B_I)$ along the last column, which involves
the determinants of the matrices $B_{I,i}$ which are the restrictions
of $B$ to the last $d-r$ lines except $i$, and to the column indexed
by the complement $I^c$ of $I$, except $d$. We thus get
\bee\label{P-22}
\ga_r(A,B)=Z+\sum_{i=1}^{d-r}(-1)^i\big(\Up_d^{r+i}+\Te^{r+i}_d\big)
Z'_i,\quad
Z'_i=\sum_{I\in\ja_r,\,d\notin I}\ep_I\det(A_I)\,\det(B_{I,i}),
\eee
where $Z$ and all $Z'_i$ and $\Te_d^{r+i}$ are independent of the vector
$\Up_d$. Since this random vector $\Up_d$ has a density, it follows that
the variable $\ga_r(A,B)$ also has a density, provided $Z'_i\neq0$ a.s.
for at least one value of $i$.

At this stage, we observe that $Z'_i$ has exactly the same structure
as $\ga_r(A,B)$, except that the dimension of each $B_{I,i}$ is
$(d-r-1)\times(d-r-1)$ instead of $(d-r)\times(d-r)$, and that the last
column of the original problem has totally disappeared. We can repeat the
argument, to obtain that $Z'_i$ has a density and is thus a.s.
non-vanishing, as soon as some similar quantity (where the last two
columns of the original problem no longer show up) is a.s. non-vanishing.
Then, after an obvious induction, we deduce that $\ga_r(A,B)$ has a
density as soon as $\det(A_I)\neq0$ a.s. for $I=\{1,\cdots,r\}$.

However, since the entries of this last $A_I$ are $\la_i\Phi^i_j$ for
$i,j=1,\cdots,r$, and all those $\la_i$ are non zero,
it is well known (and also a simple consequence of the previous proof,
in which we develop $\det(A_I)$ according to its last column and perform
the same induction procedure) that $\det(A_I)$ has a density. This indeed
shows us that $\ga_r(A,B)$ has a density, hence $\E(\ga_r(A,B)^2)>0$ and
the proof of the first part of (\ref{1-9}) is complete. \qed
\vsc

\nib 2) The results about $\Ga'_r(\uu)$ and $\Ga''_r(\uu)$. \rm When
the rank of $\al$ is smaller than $r$, we have seen that, with the previous
notation, $\ga_r(A,B)=0$, hence also $\Ga'_r(\uu)=\Ga''_r(\uu)=0$.

Next, we turn to the case when the rank of $\al$ is $r$. Exactly as in the
previous proof, it suffices to show the result when $\uu=(\La,J,\ga,a)$.
Recalling that $\BF_r(\uu,1)$ and $\BF_r(\uu,2)$ have the same law, we
have $\BE((\BF_r(\uu,1)-\BF_r(\uu,2))^2)=
2(\Ga_r'(\uu)-\Ga_r''(\uu))$ and thus the second part of (\ref{1-9})
holds unless $\BF_r(\uu,1)=\BF_r(\uu,2)$ a.s.

With the previous notation, we have $\BF_r(\uu,1)=\ga_r(A,B)^2$, and also
$\BF_r(\uu,2)=\ga_r(\BA,\BB)^2$, where $\BA$ and $\BB$
are given again by (\ref{P-21}) with the same $\la_j$'s, and random
vectors $(\BPhi_j,\BUp_j,\BTe_j)$ having globally the same distribution
as $(\Phi_j,\Up_j,\Te_j)$ (which may of course be defined for $j>d$):
these two families of vector are not independent, and we have in fact
\bee\label{P-23}
\BPhi_j=\frac1{\sqrt{2}}\,(\Phi_{2j-1}+\Phi_{2j}),\quad
\BUp_j=\frac1{\sqrt{2}}\,(\Up_{2j-1}+\Up_{2j}),
\eee
plus a more complicated relation relating $\BTe_j$ with the $\Te_{j'}$
for $j'\leq 2j$ and the vector $a$. What we need to prove is then
$\PP(|\ga_r(A,B)|\neq|\ga_r(\BA,\BB)|)>0$.

We have (\ref{P-22}), and also, by the same argument,
$$\ga_r(\BA,\BB)=\BZ+\sum_{i=1}^{d-r}(-1)^i\Big(
\frac{\Up_{2d+1}^{r+i}+\Up_{2d}^{r+i}}{\sqrt{2}}
+\Te^{r+i}_d\Big)\BZ'_i,\quad
\BZ'_i=\sum_{I\in\ja_r,\,d\notin I} \ep_I\det(\BA_I)\,\det(\BB_{I,i}),$$
where we have also used the second part of (\ref{P-23}). Here, the vector
$\Up_{2d}$ has a density and is independent of all other terms
showing in the above expression, and also independent of $\ga_r(A,B)$.
Therefore, $|\ga_r(A,B)|\neq|\ga_r(\BA,\BB)|$ almost surely on the set
$\{\BZ'_i\neq0\}$. Now, $\BZ'_i$ is the same as $Z'_i$, upon replacing
$(A,B)$ by $(\BA,\BB)$ everywhere, hence the previous proof shows that
indeed $\BZ'_i\neq0$ a.s. This shows that in fact
$\PP(|\ga_r(A,B)|\neq|\ga_r(\BA,\BB)|)=1$, thus ending the proof of the
second part of (\ref{1-9}). \qed

\subsection{Some stochastic calculus preliminaries.} \label{sec5-3}
We assume (H) and, by localization (see e.g. section 3 in \cite{BGJPS}), we may also assume that all processes
$X_t$, $\si_t$, $b_t$, $a_t$, $v_t$, $a'_t$, $v'_t$, $a''_t$, $v''_t$
are uniformly bounded in $(\omega ,t)$. The constants are always written as $K$, or $K_p$
if we want to stress the dependency on an additional parameter $p$,
and never depend on $t,i,n,j$. For any process $Y$,
we use the following simplifying notation:
\bee\label{P-109}
\f^n_i=\f_{2id\De_n},\qquad Y^n_i=Y_{2id\De_n}.
\eee
For all $p,t,s>0$, we have by Burkholder-Gundy inequality
\bee\label{P-111}
\E\big(\sup_{u\in[0,s]}\,|Y_{t+u}-Y_t|^p\mid\f_t\big)\leq K_qs^{q/2}\quad
\text{if}~Y=X,\si,b,v.
\eee
We set
\bee\label{P-120}
\eta_{t,s}=\sup_{u\in[0,s],~Y=a,v',v''}\,\|Y_{t+u}-Y_t\|^2,\qquad
\eta^n_i=\sqrt{\E(\eta_{2id\De_n,2d\De_n}\mid\f_i^n)}.
\eee

\begin{lem}\label{LP-1} For all $t>0$ we have
$\De_n\E\big(\sum_{i=0}^{[t/2d\De_n]-1}\eta^n_i\big)\to0$.
\end{lem}

\nib Proof. \rm It suffices to prove the result separately when $Y=a$ or
$Y=v'$ or $Y=v''$. Set $\gamma^n_t=
\sup_{s\in(0,4d\Delta_n]}\,\|Y_{t+s}-Y_{t}\|^2$,
so $\mathbb{E}((\eta^{n}_i)^2)$ is smaller than $\mathbb{E}(\gamma^n_0)$
when $i=0$ and than $\frac1{2d\Delta_n}\int_{2(i-1)d\Delta_n}
^{2id\Delta_n}\mathbb{E}(\gamma^n_s)\,ds$ when $i\geq1$. Hence
by the Cauchy-Schwarz inequality,
$$\Delta_n\mathbb{E}\Big(\sum_{i=1}^{[t/2d\Delta_n]-1}\eta^n_i\Big)
\leq \frac{\sqrt{t}}{\sqrt{2d}}\,\Big(\mathbb{E}\Big(
\Delta_n\sum_{i=0}^{[t/2d\Delta_n]-1}(\eta^n_i)^2\Big)\Big)^{1/2}
\leq K_t\,\Big(\mathbb{E}\Big(\gamma^n_0+
\int_0^t\gamma^n_s\,ds\Big)\Big)^{1/2}.$$
We have $\gamma^n_s\leq K$, whereas the c\`adl\`ag property of $Y$ yields
that $\gamma^n_s(\omega)\to0$ for all $\omega$, and all $s$ except for
countably many strictly positive values (depending on $\omega$). Then,
the claim follows by the dominated convergence theorem.\qed \newline \newline
The proof of Theorem \ref{T1} is based on a decomposition of the increments
$Z^n_{(2id+\ka j)\De_n}-Z^n_{(2id+\ka(j-1))\De_n}$. In order to understand
better this decomposition, we first deduce from (\ref{1-1}) that, for any
$z\leq t\leq s$, and with vector notation,
$$\int_t^sb_u\,du=\al_1+\al_2+\al_3+\al_4,\qquad
\int_t^s\si_u\,dW_u=\al_5+\al_6+\al_7+\al_8+\al_9+\al_{10}+\al_{11},$$
where
$$\begin{array}{l}
\al_1=b_z(s-t),\quad\al_2=\int_t^s\Big(\int_z^ua'_w\,dw\Big)du
\quad\al_3=v'_z\int_t^s(W_u-W_z)\,du,\\
\al_4=\int_t^s\Big(\int_z^u(v'_w-v'_z)\,dW_w\Big)\,du\\
\al_5=\si_z(W_s-W_t),\quad\al_6=a_z\int_t^s(u-z)dW_u,\quad
\al_7=\int_t^s\Big(\int_z^u(a_w-a_z)\,dw\Big)\,dW_u\\
\al_8=v_z\int_t^s(W_u-W_z)\,dW_u,\quad
\al_9=\int_t^s\Big(\int_z^u\Big(\int_z^wa''_r\,dr
\Big)\,dW_w\Big)\,dW_u\\
\al_{10}=v''_z\int_t^s\Big(\int_z^u(W_w-W_z)\,dW_w\Big)\,dW_u,\quad
\al_{11}=\int_t^s\Big(\int_z^u
\Big(\int_z^w(v''_r-v''_z)\,dW_r\Big)\,dW_w\Big)\,dW_u.
\end{array}$$
A repeated use of the Burkholder-Gundy and H\"older inequalities shows
that, in view of our assumptions on the various coefficients,
and for any $p\geq1$:
$$\E(|\al_j|^p\mid\f_z)\leq\left\{\begin{array}{ll}
K_p(s-z)^{p/2}&\text{if}~~j=5\\
K_p(s-z)^p&\text{if}~~j=1,8 \\
K_p(s-z)^{3p/2}&\text{if}~~j=3,6,10\\
K_p(s-z)^{2p}&\text{if}~~j=2,9\\
K_p(s-z)^{3p/2}\,\E(\eta_{z,s-z}^{p}\mid\f_z)&\text{if}~~j=4,7,11.
\end{array}\right.$$

We can then apply the previous decomposition with $z=2id\De_n$,
$t=(2id+\ka(j-1))\De_n$ and $s=(2id+\ka j)\De_n$, and add the increment of
the process $X'$, to obtain
\bee\label{P-110}
\frac{Z^{n,\ka}_{(2id+j\ka)\De_n}-Z^{n,\ka}_{(2id+(j-1)\ka)\De_n}}
{\sqrt{\ka\De_n}}=\al^{n,\ka}_{i,j}+\sqrt{\ka\De_n}\,\be^{n,\ka}_{i,j}
+\ka\De_n\,\ga^{n,\ka}_{i,j}+\De_n\,\de^{n,\ka}_{i,j}
\eee
for $\ka=1,2$, and where (explicitly writing the components)
$$\!\!\begin{array}{l}
\al^{n,\ka,l}_{i,j}=\frac1{\sqrt{\ka\De_n}}\,\sum\limits_{m=1}^q\si_i^{n,lm}
(W^m_{(2id+\ka j)\De_n}-W^m_{(2id+\ka(j-1))\De_n})\\
\be^{n,\ka,l}_{i,j}=b^{n,l}_i+\frac1{\ka\De_n}\,
\sum\limits_{m,k=1}^qv^{n,lmk}_i
\int_{(2id+\ka(j-1))\De_n}^{(2id+\ka j)\De_n}(W_s^k-W^{k}_{2id\De_n})
\,dW^m_s\\
\qquad+\frac1{\sqrt{\ka\De_n}}\,\sum\limits_{m=1}^d\Wsi^{lm}
(W'^m_{(2id+\ka j)\De_n}-W'^m_{(2id+\ka(j-1))\De_n})\\
\ga^{n,\ka,l}_{i,j}=\frac1{(\ka\De_n)^{3/2}}\,
\Big(\sum\limits_{m=1}^qa^{n,lm}_i
\int_{(2id+\ka(j-1))\De_n}^{(2id+\ka j)\De_n}
\big(s-2id\De_n\big)\,dW^m_s\\
\qquad+\sum\limits_{m=1}^qv'^{n,lm}_i\,
\int_{(2id+\ka(j-1))\De_n}^{(2id+\ka j)\De_n}(W^m_s-W^{m}_{2id\De_n})\,ds\\
\qquad+\sum\limits_{m,l,k=1}^qv''^{n,mlk}_i
\int_{(2id+\ka(j-1))\De_n}^{(2id+\ka j)\De_n}
\big(\int_{(2id+\ka(j-1))\De_n}^s
(W^k_u-W^{k}_{2id\De_n})\,dW^l_u\big)dW^m_s\Big)
\end{array}$$
and $\de^{n,\ka}_{i,j}$ is a remainder term, and for $p\geq1$
we have the estimates when $j\leq2d$ if $\ka=1$ and $j\leq d$ when $\ka=2$
(recalling $\eta_{t,s}\leq K$):
\bee\label{P-113}
\begin{array}{l}
\E\big(\|\al^{n,\ka}_{i,j}\|^p+\|\be^{n,\ka}_{i,j}\|^p+
\|\ga^{n,\ka}_{i,j}\|^p\mid\f_i^n\big)\leq K_p\\[1.5 ex]
\E\big(\|\de^{n,\ka}_{i,j}\|^p\mid\f_i^n\big)\leq K_p
\,(\De_n^{p/2}+(\eta^n_i)^{2\wedge p}\big)\leq K_p.
\end{array}
\eee

We end these preliminaries with a lemma which compares $S^{n,\ka}$
for $\ka=1,2$ with the following processes:
\bee\label{P-115}
\begin{array}{c}
S(r)^{n,\ka}_t=
2d\De_n\sum\limits_{i=0}^{[t/2d\De_n]-1}\ga_r(A^{n,\ka}_i,B^{n,\ka}_i)^2,
\quad\text{where}\\
A^{n,\ka}_i=\text{\rm mat}(\al^{n,\ka}_{i,1},\cdots,\al^{n,\ka}_{i,d}),\qquad
B^{n,\ka}_i=\text{\rm mat}(\be^{n,\ka}_{i,1},\cdots,\be^{n,\ka}_{i,d}).
\end{array}
\eee
It also compares $V^{n,\ka,\ka'}$ of (\ref{2-26}) with
\bee\label{P-117}
V(r)^{n,\ka,\ka'}_t=
4d^2\De_n\sum\limits_{i=0}^{[t/2d\De_n]-1}\ga_r(A^{n,\ka}_i,B^{n,\ka}_i)^2
\,\ga_r(A^{n,\ka'}_i,B^{n,\ka'}_i)^2.
\eee

\begin{lem}\label{LP-15} If $r_t(\om)\leq r$ identically for some
$r\in\{0,\cdots,d\}$, we have for $\ka,\ka'=1,2$:
\bee\label{P-118}
\rdnn\,\Big(\frac1{(\ka\De_n)^{d-r}}\,S^{n,\ka}-S(r)^{n,\ka}
\Big)~\toucp~0
\eee
and
\bee\label{P-119}
\frac1{(\ka\ka'\,\De_n^2)^{d-r}}\,V^{n,\ka,l}-V(r)^{n,\ka,\ka'}
~\toucp~0
\eee
\end{lem}

\nib Proof. \rm We denote by $\xi^{n,\ka}_i$ the $i$th summand in the
definition \eqref{1-7} of $S^{n,\ka}_t$. Besides the matrices in
(\ref{P-115}), we also define
$$C^{n,\ka}_i=\text{\rm mat}(\ga^{n,\ka}_{i,1},\cdots,\ga^{n,\ka}_{i,d}),
\qquad
D^{n,\ka}_i=\text{\rm mat}(\de^{n,\ka}_{i,1},\cdots,\de^{n,\ka}_{i,d}).$$

We start with (\ref{P-118}). Applying (\ref{P-105})
with $h=\sqrt{\ka\De_n}$, the fact that each $A^{n,\ka}_i$ has at most rank
$r$ (because $r_t\leq r$), and the estimates (\ref{P-113})
plus the Cauchy-Schwarz inequality, we obtain
$$\begin{array}{l}
\frac1{(\ka\De_n)^{d-r}}\,\xi_i^{n,\ka}=\ga_r(A^{n,\ka}_i,B^{n,\ka}_i)^2+
2\sqrt{\ka\De_n}\,\ze^{n,\ka}_i+\Wze^{n,\ka}_i,\quad\text{where}\\[1.5 ex]
\qquad\ze^{n,\ka}_i=\ga_r(A^{n,\ka}_i,B^{n,\ka}_i)
(\ga_{r-1}(A^{n,\ka}_i,B^{n,\ka}_i)
+\ga'(A^{n,\ka}_i,B^{n,\ka}_i,C^{n,\ka}_i))\\[1.5 ex]
\qquad\E(|\Wze^{n,\ka}_i|)\leq K\De_n+K\rdn\,\E(\eta^n_i).
\end{array}$$
In view of Lemma \ref{LP-1}, $\rdn\,\sum_{i=0}^{[t/2d\De_n]-1}\Wze^{n,\ka}_i
\toucp0$. Since
$S^{n,\ka}_t=2d\De_n\sum_{i=0}^{[t/2d\De_n]-1}\xi^{n,\ka}_i$, it
remains to prove that $\De_n\,\sum_{i=0}^{[t/2d\De_n]-1}\ze^{n,\ka}_i\toucp0$.

For this purpose we use the decomposition
$\ze^{n,\ka}_i=\ze'^{n,\ka}_i+\ze''^{n,\ka}_i$, where
$\ze'^{n,\ka}_i=\E(\ze^{n,\ka}_i\mid\f_i^n)$. By Doob's inequality,
(\ref{P-113}) and the fact that $\ze''^n_i$ is $\f_{i+1}^n$-measurable, we
have
$$\E\Big(\sup_{s\leq t}\Big(\sum_{i=0}^{[s/2d\De_n]-1}
\ze''^{n,\ka}_i\Big)^2\Big)\leq
2^{d+1}\,\E\Big(\sum_{i=0}^{[t/2d\De_n]-1}|\ze^n_i|^2\Big)\leq
\frac{Kt}{\De_n}.$$
Thus $\De_n\,\sum_{i=0}^{[t/2d\De_n]-1}\ze''^{n,\ka}_i\toucp0$, and the
result will hold if we can prove that $\ze'^{n,\ka}_i=0$.
We even prove the stronger statement that
$\E(\ze^{n,\ka}_i\mid\g^{W'}\vee\f_i^n)=0$, where $\g^{W'}$ is the
$\si$-field generated by the whole process $W'$, and this is implied by
\bee\label{P-116}
\begin{array}{l}
\III\in\ia_{(r,d-r)},~~
\III'\in\ia_{(r-1,d-r+1)},~~\III''\in\ia_{(r,d-r-1,1)}~~\Longrightarrow\\[1.5 ex]
\qquad \E\big(\det(G^\III_{A^{n,\ka}_i,B^{n,\ka}_i})\,
\det(G^{\III'}_{A^{n,\ka}_i,B^{n,\ka}_i})\mid
\g^{W'}\vee\f_i^n\big)=0\\ [1.5 ex]
\qquad \E\big(\det(G^\III_{A^{n,\ka}_i,B^{n,\ka}_i})\,
\det(G^{\III''}_{A^{n,\ka}_i,B^{n,\ka}_i,C^{n,\ka}_i})\mid
\g^{W'}\vee\f_i^n\big)=0.
\end{array}
\eee
The variables $\al^{n,\ka,l}_{i,j}$, $\be^{n,\ka,l}_{i,j}$ and
$\ga^{n,\ka,l}_{i,j}$ have
the form $\Phi(\om,(W(\om)_{2id\De_n+t}-W(\om)_{2id\De_n})_{t\geq0})$,
with $\Phi$ a $(\g^{W'}\vee\f_i^n)\otimes\ca^d$-measurable
function on $\Om\times C(\R_+,\R^d)$, where $C(\R_+,\R^d)$ is the set of all
continuous $\R^d$-valued functions on $\R_+$ and $\ca^d$ is its Borel
$\si$-field for the local uniform topology. When $\Phi=
\al^{n,\ka,l}_{i,j}$ or $\Phi=\ga^{n,\ka,l}_{i,j}$, the map
$x\mapsto\Phi(\om,x)$ is odd, in the sense that $\Phi(\om,-x)=
\Phi(\om,x)$, and it is even when $\Phi=\be^{n,\ka,l}_{i,j}$.

In (\ref{P-116}), the three variables
$\det(G^\III_{A^{n,\ka}_i,B^{n,\ka}_i})$,
$\det(G^{\III'}_{A^{n,\ka}_i,B^{n,\ka}_i})$,
$\det(G^{\III''}_{A^{n,\ka}_i,B^{n,\ka}_i,C^{n,\ka}_i})$
are associated with three functions $\Psi$, $\Psi'$, $\Psi''$ of the same
type. What precedes yields that $\Psi$ is even (resp. odd)
if $r$ is even (resp. odd), and both $\Psi'$ and $\Psi''$ are even (resp. odd)
if $r$ is odd (resp. even). Consequently, the products $\Psi\Psi'$
and $\Psi\Psi''$ are odd in all cases. Since the
$\g^{W'}\vee\f_i^n$-conditional law of
$(W_{2id\De_n+t}-W_{2id\De_n})_{t\geq0}$ is invariant by the map
$x\mapsto-x$ on $C(\R_+,\R^d)$, we deduce (\ref{P-116}), hence (\ref{P-118})
holds.

Finally, we turn to (\ref{P-119}). Let $\te^n_i$ be the $i$th summand
in the right side of \eqref{2-26}, for $\ka,\ka'$ fixed. We can apply
(\ref{P-106}) with $h=\sqrt{\ka\De_n}$ and $h=\sqrt{\ka'\De_n}$, and
\eqref{P-110} and (\ref{P-113}) again, to get
$$\E\Big(\Big|\frac1{(\ka\ka'\De_n^2)^{d-r}}\,\te^n_i
-\ga_r(A^{n,\ka}_i,B^{n,\ka}_i)^2\,
\ga_r(A^{n,\ka'}_i,B^{n,\ka'}_i)^2\Big|\Big)\leq K\rdn.$$
(\ref{P-119}) follows, and the proof is complete.\qed

\subsection{Proof of Theorem \ref{T1}, Corollary \ref{C12},
and Proposition \ref{P3}.} \label{sec5-4}

\nib 1) \rm Observe that,
with the notation (\ref{1-31}) and \eqref{1-32}, and upon taking
\bee\label{P-5}
\uu^n_i=(\si_i^n,\Wsi,v_i^n,b_i^n),\quad
\BW_t=\frac{W_{(2id+t)\De_n}-W_{2id\De_n}}{\rdn},\quad
\BW'_t=\frac{W'_{(2id+t)\De_n}-W'_{2id\De_n}}{\rdn},
\eee
we have $\ga_r(A^{n,\ka}_i,B^{n,\ka}_i)^2=\BF(\uu^n_i,\ka)$. We consider the
two-dimensional variables $\xi^n_i$ with components
\begin{align} \label{xi}
\xi^{n,\ka}_i=2d\,\rdn\,\big(\ga_r(A^{n,\ka}_i,B^{n,\ka}_i)^2
-\Ga_r(u^n_i)\big), \qquad \ka=1,2.
\end{align}
Since $\uu^n_i$ is $\f^n_i$-measurable, whereas the processes
$\BW$ and $\BW'$ above are independent of $\f^n_i$, we deduce from
\eqref{1-8}, and from \eqref{P-113} for the estimate below, that
\bee\label{P-7}
\begin{array}{l}
\E(\xi^{n,\ka}_i\mid\f_i^n)=0,\qquad \E(\|\xi^n_i\|^4\mid\f_i^n)\leq K\De_n^2\\[1.5 ex]
\E(\xi^{n,\ka}_i\,\xi^{n,\ka'}_i\mid\f_i^n)=\left\{\begin{array}{ll}
4d^2\De_n\,\Ga'_r(\uu_i^n)&\text{if}~\ka=\ka'\\
4d^2\De_n\,\Ga''_r(\uu_i^n)&\text{if}~\ka\neq\ka'\end{array}\right.
\end{array}
\eee
\qed
\vsq

\nib 2) \rm By \eqref{P-118}, for Theorem \ref{T1}
it is enough to prove the stable convergence $U'(r)^n\tolls \ua(r)$, where
$U'(r)^n$ is the two-dimensional process with components
$U'(r)^{n,\ka}=\rdnn\,(S(r)^{n,\ka}-S(r))$ and the quantity $S(r)^{n,\ka}$ is
defined in \eqref{P-115}. We have $U'(r)^n=Y^n+Y'^n$, where
$$\begin{array}{l}
Y^n_t=\sum\limits_{i=0}^{[t/2d\De_n]-1}\xi^n_i\\
Y'^n_t=\rdnn\,\Big(2d\De_n\sum\limits_{i=0}^{[t/2d\De_n]-1}
\Ga_r(\si^n_i,\Wsi,v^n_i,b^n_i)-
\int_0^t\Ga_r(\si_s,\Wsi,v_s,b_s)\,ds\Big),
\end{array}$$
and $\xi^n_i$ is given in \eqref{xi}.
Since the three processes $\si,v,b$ are It\^o semimartingales, whereas
$\Ga_r$ is a $C^\infty$ function, it is well known that $Y'^n\toucp0$ (see e.g. section 8 in \cite{BGJPS}).
We are thus left to prove that
\bee\label{P-8}
Y^n~\tolls~\ua(r).
\eee
By virtue of the first two parts of \eqref{P-7}, a standard CLT
for triangular arrays of martingale (see \cite[Theorem IX.7.28]{JS}) increments shows that, for
(\ref{P-8}) to hold, it suffices to show the next two properties:
\bee\label{P-10}
\sum_{i=0}^{[t/2d\De_n]-1}\E(\xi^{n,\ka}_i\,\xi^{n,\ka'}_i\mid\f_i^n)
~\toop~V(r)^{\ka\ka'}_t
\eee
\bee\label{P-11}
\sum_{i=0}^{[t/2d\De_n]-1}\E\big(\xi^n_i\,(N_{(2(i+1)d\De_n}-
N_{2id\De_n})\mid\f_i^n\big)~\toop~0
\eee
for all $t>0$ and for any bounded martingale $N$ orthogonal to $(W,W')$ and
also for $N=W^m$ or $N=W'^m$ for any $m$.

The last part of (\ref{P-7}) and the c\`adl\`ag property of $\si,v,b$,
plus the fact that $\Ga'_r$ and $\Ga''_r$ are polynomials, immediately gives
us \eqref{P-10} by Riemann integration.

The proof of \eqref{P-11} is also standard: By construction,
$\xi^n_i$ is a two-dimensional variable of the
form $\Phi(\om,(W_{2id\De_n+t}(\om)-W_{2id\De_n}(\om))_{t\geq0},
(W'_{2id\De_n+t}(\om)-W'_{2id\De_n}(\om))_{t\geq0})$ similar to the functions
occurring in Lemma \ref{LP-15}, and since in the definition of $\xi^n_i$
one takes squared determinants, all these functions $\Phi$ are globally
even in the sense
that $\Phi(\om,x,y)=\Phi(\om,-x,-y)$ for any two $d$-dimensional functions
$x,y$ on $\R_+$. So, on the one hand, after multiplying the function $\Phi$
corresponding to $\xi^n_i$ by $x^m(2d\De_n)$ or
$y^m(2d\De_n)$, one gets an odd function, and (\ref{P-11})
when $N=W^m$ or $N=W'^m$ follows. On the other hand, by the representation
theorem one can write $\xi^n_i$ as the sum of two integrals over
$(2id\De_n,2(i+1)d\De_n]$ with respect to $W$ and $W'$, for suitable
predictable integrands; thus when $N$ is orthogonal to $W$ and $W'$,
the increment $N_{(2(i+1)d\De_n}-N_{2id\De_n}$ has $\f_i^n$-conditional
correlation $0$ with both those integrals, thus yielding (\ref{P-11}) again.

Therefore, the proof of Theorem \ref{T1} is complete. \qed
\vsq

\nib 3) \rm A simple calculation shows that
\[
\wR(n,T)-r=\frac1{\log2}\,
\log\frac{1+\rdn\,U(r)^{n,1}_T/S(r)_T}{1+\rdn\,U(r)^{n,2}_T/S(r)_T} \quad \text{if} \quad
S(r)_T>0,
\]
hence on the set $\Om^r_T$. Since the sequence $U(r)^n_T$
is tight, it follows from a Taylor expansion that
\bee\label{P-13}
\rdnn\,(\wR(n,T)-r)-\frac1{S(r)_T\,\log2}\,\big(
U(r)^{n,1}_T-U(r)^{n,2}_T\big)~\toop~0
\eee
on $\Om_T^r$ again. Then Corollary \ref{C12} follows from Theorem \ref{T1},
upon observing that the $\f$-conditional variance of $U(r)^1_T-U(r)^2_T$ is the
numerator of the right side if \eqref{2-7}. \qed
\vsq

\nib 4) \rm Now we turn to the proof of \eqref{2-29}, and by (\ref{P-119})
it suffices to prove the convergence of $V(r)^{n,\ka,\ka'}$.
We suppose that $\ka=\ka'$, the proof in the case $\ka\neq\ka'$ being
analogous. We set
$$\eta^n_i=\ga_r(A_i^{n,\ka},B_i^{n,\ka})^4,\qquad
\eta'^n_i=\E(\eta^n_i\mid\f^n_i),\qquad\eta''^n_i=\eta^n_i-\eta'^n_i.$$
As for \eqref{P-7}, we deduce from \eqref{1-8} and \eqref{P-113} that
$$\eta'^n_i=2d\big(\Ga_r'(\uu^n_i)-\Ga_r(\uu^n_i)^2\big),\qquad
\E(|\eta''^n_i|^2)\leq K.$$
On the one hand, the same argument as for proving \eqref{P-10} shows
that $4d^2\De_n\sum_{i=0}^{[t/2d\De_n]-1}\eta'^n_i$ converges in the u.c.p.
sense to the right side of \eqref{2-29} (for $\ka=\ka'$). On the other hand,
since $\eta''^n_i$ is a martingale increment relative to the filtration
$(\f^n_i)_{i\geq0}$, we deduce from Doob's inequality that
$4d^2\De_n\sum_{i=0}^{[t/2d\De_n]-1}\eta''^n_i\toucp0$. We then deduce
\eqref{2-29}. \qed
\vsq

\nib 5) \rm Finally, for \eqref{2-30} it is enough to show the convergence
in probability in restriction to each set $\Om^r_T$, for $r=0,\cdots,d$.
For this we use the following convergence properties, which readily follow
from \eqref{2-1}, \eqref{2-6}, in restriction to the set $\Om^r_T$:
$$\frac1{\De_n^{d-r}}\,S^{n,1}\toop S(r)_T>0,\qquad
\wR(n,T)\toop r,$$
together with \eqref{2-29} applied at time $T$, which also holds on $\Om^r_T$.
Then \eqref{2-30} follows after a (slightly tedious) calculation, in view of
the form \eqref{2-7} of $V(T)$ on $\Om^r_T$: the proof of
Proposition \ref{P3} is complete. \qed

\subsection{Proof of Theorems \ref{T4} and  \ref{T5}.} \label{sec5-5}

We begin with a lemma. Its setting apparently extends the setting of the
theorem to be proved, but this will be useful for the proof itself. The
extension concerns
the fact that we replace the non-random terminal time $T$ by a stopping
time, still denoted by $T$, which is positive and bounded. In this case, the
notation \eqref{1-2} still makes sense, as well as
$A(p)^n_T$ and $a(n,T)$ and $\BV_T^{n,\ka\ka'}$, as given by \eqref{2-69}
and \eqref{2-50}.

\begin{lem}\label{L-7} Assume (H) and $r_t=r$ for all $t\leq T$ with $T$ a
positive finite stopping time and $r\in\{0,\cdots,d\}$. Then for all $p>0$,
$\ka,\ka'\in\{1,2\}$ and $\Te^{r,\ka\ka'}_s$ as in \eqref{2-29} we have
\bee\label{P-90}
A(p)^n_T~\toop~T\,r^p
\eee
\bee\label{P-96}
\frac1{(\ka\ka'\De_n^2)^{d-r}}\,\BV^{n,\ka\ka'}_T~\toop~
2d\int_0^T\Big(\frac1{\Ga_r(\si_s,\Wsi,v_s,b_s)}-\frac T{S(r)_T}\Big)^2
\,\Te^{r,\ka\ka'}_s\,ds.
\eee
Moreover, if $r\geq1$, the following stable convergence in law holds,
where $\ua(r)$ is defined in Theorem \ref{T1}:
\bee\label{P-91}
\begin{array}{l}
\Big(U(r)^n_T,\rdnn\,(A(p)^n_T-a(n,T)\,r^p)\Big)~\tols\\
\hskip4cm\Big(\ua(r)_T,\frac{pr^{p-1}}{\log 2}\int_0^T
\frac1{\Ga_r(\si_s,\Wsi,v_s,b_s)}\,(d\ua(r)^1_s-d\ua(r)^2_s)\Big).
\end{array}
\eee
\end{lem}

\nib Proof. 1) \rm Let $\ga_t=\Ga_r(\si_t,\Wsi,v_t,b_t)$, which is a
continuous process, positive on $[0,T]$ by Lemma \ref{L-0}. Thus
$T_m=m\wedge T\wedge\inf(t:\ga_t<1/m)$ satisfies $\PP(T_m=T)\to1$ as
$m\to\infty$ and, if any one of the claimed convergence holds for each
$T_m$ (instead of $T$), it also holds for $T$. In other words, we can assume
$T\leq A$ and $1/\ga_t\leq A$ for some constant
$A$ and all $t\in[0,T]$. Moreover, $\Ga_r$ is a polynomial, so
the process $\ga_t$ is a continuous It\^o semimartingale, and by localization
again one can assume that for some other constant $A'$,
\bee\label{P-80}
\E(|\ga_{t+s}-\ga_t|^2)\leq A'\,s.
\eee
The sequence $U(r)^n$ converges in law toward a continuous process, so the
moduli of continuity $\rho(n,x)=\sup(\|U(r)^n_{t+s}-U(r)_t^n\|:~
t\leq A',~|s|\leq x)$ satisfy $\lim_{x\downarrow0}\,\limsup_n\,\PP(\rho(n,x)>1)
=0$, and thus with the simplifying notation $w_n=2dk_n\De_n$ we have
\bee\label{P-83}
\PP(\Om_n)\to1,\qquad\text{where}~~
\Om_n=\{\|U(r)^n_{t+s}-U(r)^n_t\|\leq 1~\forall\,t\leq A',~s\leq w_n\big\}.
\eee
\vsd

\nib 2) \rm Observe that
$$\begin{array}{c}
\wR^n_i=r+\frac1{\log2}\,\log\frac{\ze^n_i
+\rdn\,\eta^{n,1}_i}{\ze^n_i+\rdn\,\eta^{n,2}_i},\quad\text{where}\\[1.5 ex]
\ze^n_i=S(r)_{2id\De_n+w_n}-S(r)_{2id\De_n},\quad
\eta^{n,k}_i=U(r)^{n,k}_{2id\De_n+w_n}-U(r)^{n,k}_{2id\De_n}.
\end{array}$$
Recalling $1/\ga_t\leq A$, and since
$\ze^n_i=w_n(\ga_{2id\De_n}+\rho^n_i)$, where
$\E(|\rho^n_i|^2)\leq A''w_n$ by \eqref{P-80}, one has
\bee\label{P-99}
0\leq \frac{w_n}{\ze^n_i}\leq  A,\qquad
\E\Big(\Big|\frac{w_n}{\ze^n_i}-\frac1{\ga_{2id\De_n}}\Big|^2\Big)
\leq A^2A''w_n.
\eee
Moreover, take $\al\in(0,1/2)$ such that
$\frac{2|\log(1-\al)|}{\log 2}\leq1/2$. For $n$ large enough we have
$A\,\rdn\,/w_n\leq \al$ because $k_n\De_n^{3/4}\to\infty$. In this case,
in restriction to the set $\Om_n$,
for all $i\leq[T/2d\De_n]-k_n-1$ we have with a constant $K$ (varying from
place to place below):
\bee\label{P-93}
\begin{array}{ll}
\Big|\frac{\rdn\,\eta^{n,k}_i}{\ze^n_i}\Big|\leq \frac{A\,\rdn}{w_n}
\leq\frac12,\qquad&
\Big|\wR^n_i-r-\frac{\rdn}{\log2}\,\frac{\eta^{n,1}_i-\eta^{n,2}_i}
{\ze^n_i}\Big|\leq K\,\frac{\De_n}{w_n^2}\\[1.5 ex]
r\geq1~\Rightarrow~\Big|\frac{\wR^n_i}r-1\Big|
\leq\frac{K\,\rdn}{w_n}\bigwedge\frac12,\quad&
r=0~\Rightarrow~|\wR^n_i| \leq\frac{K\,\rdn}{w_n}\bigwedge\frac12.
\end{array}
\eee

\nib 3) \rm Recalling \eqref{2-69} and $T\leq A$, when $r=0$ the last
estimate above yields
$$\E(A(p)^n_T\,1_{\Om_n})\leq KA\,\frac{\De_n^{p/2}}{w_n^p},$$
which goes to $0$ because $k_n\rdn\to\infty$. Thus in view of
\eqref{P-83} one gets \eqref{P-90} when $r=0$.
\vst

\nib 4) \rm At this stage, we start proving \eqref{P-91}, and thus
assume $r\geq1$. We observe that
$$Y_n:=\rdnn\,(A(p)^n_T-a(n,T)\,r^p)=w_n\!\!\sum_{i=0}^{[T/w_n]-2}\xi^n_i,
~~\text{where}~
\xi^n_i=\rdnn\,\big(|\wR^n_{ik_n}|^p\wedge(d+1)^p-r^p\Big).$$
\eqref{P-93} implies that for $n$ large enough, $|\wR^n_i|\leq d+1$ (recall
$r\leq d$), hence a Taylor expansion of the function $x\mapsto|r+x|^p-r^p$
imply, again for $n$ large enough:
$$\Big|\xi^n_i-\frac{pr^{p-1}}{\log2}\,\frac{\eta^{n,1}_{ik_n}
-\eta^{n,2}_{ik_n}}
{\ze^n_{ik_n}}\Big|\leq K\frac{\rdn}{w_n^2}\quad \text{on}~\Om_n~
\text{and for $i\leq[T/w_n]-2$}.$$
Upon using \eqref{P-99}, and by the Cauchy-Schwarz
inequality, it follows that
$$\E\Big(\Big|\xi^n_{i}-\frac{pr^{p-1}}{\log2}\,
\frac{\eta^{n,1}_{ik_n}-\eta^{n,2}_{ik_n}}
{w_n\,\ga_{(i-1)w_n}}\Big|\,1_{\Om_n}\Big)\leq
K\frac{\rdn}{w_n^2}+K\sqrt{w_n},$$
hence
$$\E\big(|Y_n-Y'_n|\,1_{\Om_n}\big)\to0,\quad
\text{where}~~Y'_n=\frac{pr^{p-1}}{\log2}\,\sum_{i=0}^{[T/w_n]-2}
\frac{\eta^{n,1}_{ik_n}-\eta^{n,2}_{ik_n}}{\ga_{(i-1)w_n}},$$
because $k_n\De_n^{3/4}\rightarrow \infty$.
Recall also \eqref{P-13}. Then, by virtue of \eqref{P-83}, the convergence
\eqref{P-91} will follow from
\bee\label{P-84}
\big(U(r)^n_T,Y'_n)~\tols~(\ua(r)_T,\y),\quad
\y=\frac{pr^{p-1}}{\log 2}\int_0^T
\frac1{\Ga_r(\si_s,\Wsi,v_s,b_s)}\,(d\ua(r)^1_s-d\ua(r)^2_s).
\eee
\vsd

\nib 5) \rm By Theorem VI.6.15
of \cite{JS} it follows from \eqref{P-7} and \eqref{P-10} that not only
does the sequence of processes $U(r)^n$ converge in law, but it also enjoys
the so-called P-UT property (predictable uniform tightness). By a trivial
extension of Theorem VI.6.22 in \cite{JS}, this implies that if a sequence
$H^n$ of adapted c\`adl\`ag two-dimensional processes on $\Om$
is such that the pair $(U(r)^n,H^n)\tolls(\ua(r),H)$
(functional convergence for the Skorokhod topology), the bi-dimensional
processes $\big(U(r)^n,\int_0^tH^n_{s-}\,dU(r)^n_s\big)$ converge stably in
law to $\big(\ua(r),\int_0^tH_{s-}\,d\ua(r)_s\big)$, and since $\ua(r)$ is
continuous and $T$ is $\f$-measurable, this in turn implies the
stable convergence of the variables
$\big(U(r)^n_T,\int_0^T H^n_{s-}\,dU(r)^n_s\big)\tols
\big(\ua(r)_T,\int_0^TH_{s-}\,d\ua(r)_s\big)$.
At this point, \eqref{P-84} follows, by taking the processes $H^n$ and $H$
with components
$$H^{n,1}_t=-H^{n,2}_t=\frac{pr^{p-1}}{\ga_{(i-1)w_n}\log2}
~~\text{if}~t\in((i-1)w_n,(iw_n)\wedge T],\quad H^1_t=-H^2_t=
\frac{pr^{p-1}}{\ga_t\log2}\,1_{\{t\leq T\}}.$$
(Note that the joint stable convergence $(U(r)^n,H^n)\tolls(\ua(r),H)$ holds
because $1/\ga_t$ is continuous.) This ends the proof of \eqref{P-91}.
\vst

\nib 6) \rm Since \eqref{P-91} implies \eqref{P-90} when $r\geq1$, we
are left to prove \eqref{P-96}. We fix $\ka,\ka'$.
Our first observation is that, since $U(r)^n_T\tols\ua(r)_T$ follows from
\eqref{2-1} as seen before, the proof of \eqref{2-6} carries over to the case
$T$ is a stopping time. Therefore $(\wR(n,T)-r)\,\log\De_n\toop0$ because
here $\Om^r_T=\Om$, and thus $\De_n^{\wR(n,T)-r}\toop1$. It follows that
\eqref{P-96} amounts to proving the same result for the variable
$\WV^{n,\ka\ka'}_T$ which is the same as $\BV^{n,\ka\ka'}_T$ except that
in front of the sum we substitute $\De_n^{1+2d-2\wR(n,T)}$ with
$\De_n^{1+2d-2r}$.

With $\te^n_i$ being as for Proposition \ref{P3},
the $i$th summand in the right side of \eqref{2-26}, we have
\bee\label{P-970}
\begin{array}{c}
\WV^{n,\ka\ka'}_T=\sum\limits_{j=0}^2
\Big(\frac{T\De_n^{d-r}}{S^{n,1}_T}\Big)^j
B(j)^n_T,\quad\text{where}~~
B(j)^n_T=4d^2\De_n
\sum\limits_{i=0}^{[T/2d\De_n]-k_n-1}\up(j)^n_i\te^n_i\\
\up(0)^n_i=\Big(\frac{w_n}{\ze^n_i+\rdn\,\eta^{n,1}_i}\Big)^2,
\qquad \up(1)^n_i=-2\,\frac{w_n}{\ze^n_i+\rdn\,\eta^{n,1}_i},
\qquad\up(2)^n_i=1.
\end{array}
\eee
Combining \eqref{P-99} and \eqref{P-93}, we obtain
for $i\leq[T/2d\De_n]-k_n-1$ and all $n$ large enough:
$$\E\Big(\Big|\up(0)^n_i-\frac1{(\ga_{2id\De_n})^2}\Big|^2\,1_{\Om_n}\Big)
+\E\Big(\Big|\up(1)^n_i+\frac2{\ga_{2id\De_n}}\Big|^2\,1_{\Om_n}\Big)
\leq K\Big(\frac{\De_n}{w_n^2}+w_n\Big).$$
Since by localization we may assume that the processes $\si_t,v_t,b_t$ are
bounded, we may also assume $\te^n_i\leq K$, and upon using \eqref{P-83}
once more, we then deduce that $B(j)^n_T$ as the same asymptotic behavior
as $B'(j)^n_T$ which is given by the same formula, with $\up(j)^n_i$
substituted with the following variables $\up'(j)^n_i$:
$$\up'(0)^n_i=\frac1{(\ga_{2id\De_n})^2},
\qquad \up'(1)^n_i=-2\,\frac1{\ga_{2id\De_n}},
\qquad\up'(2)^n_i=1.$$
This allows us to get, with $\Up(0)=\frac1{\ga^2}$, $\Up(1)=\frac1{\ga}$
and $\Up(2)=1$:
\bee\label{P-97}
\frac1{(\ka\ka'\De_n^2)^{d-r}}\,B(j)^n_T\toop 2d\int_0^T\Up(j)_s
\Te_s^{r,\ka\ka'}\,ds
\eee
(indeed, the case $j=2$ is \eqref{2-29}, and the other two cases follow
from a standard argument, similar to Step 5 above, but simpler because
the integrand is a c\`adl\`ag bounded process not depending on $n$, and
$\te^n_i\geq0$). Using further $T\De_n^{d-r}/S_T^{n,1}\toop T/S(r)_T$ and
recalling \eqref{P-970}, and upon expanding the square in the right side
of \eqref{P-96}, we obtain this convergence and the lemma is proved.\qed
\vsc

\nib Proof of (a) of Theorem \ref{T4}. \rm
Since $\wR(n,T)\toop R_T$ by \eqref{2-6}, it suffices to prove that
\bee\label{P-85}
A(p)^n_T~\toop~A(p)_T:=\int_0^T(r_s)^p\,ds.
\eee
The assumption implies the existence of a sequence
of stopping times $\tau_j$ increasing to infinity, such that $\tau_0=0$
and $\tau_j<\tau_{j+1}$ if $\tau_j<\infty$, and such that the process
$r_t$ takes a constant (random) value $\rho(j)$ on the time interval
$J_j=(\tau_{j-1},\tau_{j})$, with $\rho(j)\neq \rho(j+1)$ if $0<\tau_j<\infty$.
In view of the discussion preceding \eqref{1-2}, the values $r_{\tau_j}$
is necessarily smaller than or equal to $\rho(j)\wedge \rho(j+1)$, but is
irrelevant to our discussion. We also denote by $N_T$ the biggest $j$ such
that $\tau_j\leq T$.

With an empty sum being set to $0$, we have
$$A(p)^n_T=\sum_{j=1}^{N_T}Y(j)_n+Z_n,\qquad
Y(j)_n=w_n\sum_{i=[\tau_{j-1}/w_n]+1}^{[(\tau_j\wedge T)/w_n]-2}|\wR^n_i|^p,$$
and where $Z_n$ is the sum of at most $3N_T$ terms of the form
$w_n(|\wR^n_i|^p\wedge(d+1)^p)$. Since $N_T$ is finite and $w_n\to0$,
we have $Z^n_t\to0$ (pointwise), and it suffices to show that
for each $j\geq1$ we have
\bee\label{P-86}
Y(j)_n~\toop~Y(j):=\big((T\wedge\tau_{j})-(T\wedge\tau_{j-1})\big)
\rho(j)^p.
\eee
We then fix $j$. The variable $Y(j)_n$ is the process $A(p)^n$
evaluated at time $T_j=T\wedge\tau_j-\tau_{j-1}$ relative to
the underlying process $X(j)_t=X_{\tau_{j-1}+t}$, up to at most two
border terms. We thus might be tempted to apply \eqref{P-90} right away, and
indeed $X(j)$ satisfies (H) for the filtration $\f(j)_t=\f_{\tau_{j-1}+t}$,
relative to which $T_j$ is a positive bounded stopping time.
There are, however, a few problems to overcome:
\begin{enumerate}
\item The rank $r_t(X(j))$ associated with $X(j)$ is equal to $\rho(j)$
for all $t\in(0,T_j)$, but not necessarily for $t=0$, nor for $t=T_j$;
\item This rank $\rho(j)$ is random, albeit $\f(j)_0$-measurable;
\end{enumerate}
To solve these problems we fix $\ep>0$ and consider the process
$X(j,\ep)_t=X_{\tau_{j-1}+\ep+t}$, satisfying (H) for the filtration
$\f(j,\ep)_t=\f_{\tau_{j-1}+\ep+t}$, and the $\f(j,\ep)_t$-stopping
time $T(j,\ep)=T\wedge\tau_j-T\wedge\tau_{j-1}-2\ep$. The associated rank is
thus $\rho(j)$ for all $t\in[0,T(j,\ep)]$, and we will show that if
$A(p,j,\ep)$ is associated with $X(j,\ep)$ by \eqref{2-69}, we have
\bee\label{P-1001}
A(p,j,\ep)^n_{T(j,\ep)}~\toop~T(j,\ep)\,\rho(j)^p~=~
(T\wedge\tau_j-T\wedge\tau_{j-1}-2\ep)\,\rho(j)^p.
\eee
Indeed, it suffices to prove this in restriction to each set
$\Om'_r=\{\rho(j)=r\}$ satisfying $\PP(\Om'_r)>0$. If $\PP_r$ denotes
the (usual) conditional probability $\PP(\cdot\mid\Om'_r)$, the process
$X(j,\ep)$, on the space $(\Om'_r,\f\cap\Om'_r,(\f(j,\ep)_t\cap\Om'_r),\PP_r)$,
still satisfies (H) and the associated rank is now $r$ on the time
interval $[0,T(j,\ep)]$. Then Lemma \ref{L-7} yields the convergence
\eqref{P-1001} under $\PP_r$, hence under $\PP$ in restriction to each
$\Om'_r$, hence under $\PP$ on $\Om$ itself.

Finally, the difference $Y(j)_n-
A(p,j,\ep)^n_{T(j,\ep)}$ is a sum of at most $2[\ep/w_n]$ terms, each one
smaller than $w_n(d+1)^p$, so this difference is smaller than $K\ep$, as is
the difference between the two right sides of \eqref{P-86} and \eqref{P-1001}.
Hence \eqref{P-86} follows from \eqref{P-1001}, by taking first $n\to\infty$
and then $\ep\to0$. This completes the proof.\qed
\vsc

\nib Proof of (b) of Theorem \ref{T4}. \rm Exactly as in the previous
proof, it is enough to prove the result when $r_t=r\geq1$ identically, for some
non-random $r\in\{1,\cdots,d\}$. By a standard localization procedure
we can assume that $\Ga_r(\si_t,\Wsi,v_t,b_t)$, which is positive
everywhere, is bounded from below by a constant $1/A$ with $A>0$, so the assumptions
of Lemma \ref{L-7} are satisfied. Therefore, \eqref{P-13} and \eqref{P-91}
yield that, with $Y_n$ and $\y$ as in the proof of Lemma \ref{L-7} and with
$Z_n=\rdnn\,(\wR(n,T)-r)$ and $\z=\frac1{S(r)_T\,\log2}\,
(\ua(r)^1_T-\ua(r)^2_T)$, we have
\bee\label{P-81}
(Y_n,Z_n)~\tols~(\y,\z).
\eee
Then we obtain
$$\rdnn\,B(n,p,T)=Y_n+\frac{a(n,T)}{\rdn}\,\big(r^p-
\big|r+\rdn\,Z_n\big|^p\big).$$

On the one hand, $a(n,T)\to T$. On the other hand, since $Z_n$ converges
in law and $r\geq1$, we have by the mean value theorem
$$\frac{1}{\rdn}\,\big(r^p-\big|r+\rdn\,Z_n\big|^p\big)
+pr^{p-1}\,Z_n\toop0.$$
Hence \eqref{P-81} yields
$$\rdnn\,B(n,p,T)~\tols~\ba(p,T)=\y-Tpr^{p-1}\z.$$
The pair $(\y,\z)$ being $\f$-conditionally centered Gaussian, the same
is true of $\ba(p,T)$, and the form \eqref{2-72} of its conditional
variance is easily checked, by virtue of \eqref{2-5}.\qed
\vsc

\nib Proof of Theorem \ref{T5}. \rm It is easy to construct a process
$X'$ which satisfies the assumptions of (a) of Theorem \ref{T4} and
such that $X'_t=X_t$ for all $t\leq T$ on the set $\Om^{\neq}_T$. Then
on this set $B(n,p,T)$ is the same, when constructed upon $X$ or upon $X'$,
and thus it converges in probability on this set to a strictly negative
variable. On the other hand, $Z(n,p,T)$ is $B(n,p,T)$ divided by a
quantity which by construction is smaller than $\De_n^{1/4}$. Then
the convergence $Z(n,p,T)\to-\infty$ on $\Om^{\neq}_T$ is clear.

It suffices to prove \eqref{2-52} on the set
$\Om_T^=\cap\Om_T^r$ for any $r\in\{0,\cdots,d\}$ such that
$\PP(\Om_T^r)>0$. For this we can argue under the conditional probability
$\PP_r=\PP(\cdot\mid\{r_0=r\})$, or equivalently suppose that we have
in fact $r_0=r$. As above, one can construct a process $X'$ which satisfies
the assumptions of Lemma \ref{L-7} for some stopping time $T'$ which
satisfies $T'\geq T$ on the set $\Om^=_T$, and we can apply \eqref{P-96}
to $X'$ and the stopping time $T'\wedge T$. This gives us \eqref{P-96} for
$X$, in restriction to the set $\Om^=_T$.

At this point, \eqref{2-52} follows from \eqref{P-96} by exactly the
same calculations as \eqref{2-30} follows from \eqref{2-29}.

Finally, since $\BV(n,p,T)\toop \BV(p,T)$ on $\Om^=_T$, we have
$Z(n,p,T)=B(n,p,T)/\sqrt{\De_n\,\\BV(n,p,T)}$ on a set $\Om''_n$ whose
probability goes to $1$. The first part of \eqref{2-54} than follows
from \eqref{2-71} and \eqref{2-52} by delta method for stable convergence. \qed

\subsection{Proof of Corollary \ref{C3}.}
The same stopping argument as in Step 2 of the previous proof allows us
to show that, without assumptions on the rank process $r_t$, the stable
convergence in law (\ref{2-1}) holds in restriction to the set $\Om_T^r$,
as soon as we restrict our attention to the time interval $[0,T]$.

At this stage, the claim of Corollary \ref{C3} follows from \eqref{2-6}, an application
of the delta method,
\eqref{2-30} and classical properties of stable convergence in law. \qed

\subsection{Proof of Propositions \ref{P1}, \ref{P2} and \ref{P20}.}

(\ref{2-61}) is an obvious consequence of the stable convergence
(\ref{2-67}). For the alternative-consistency, it suffices to prove that
for any $r'\neq r$ we have
\bee\label{P-70}
\PP\big(\ca(\al)_T^{n,=r}\cap\Om_T^{r'}\big)\to\PP(\Om_T^{r'}).
\eee
On the set $\Om_T^{r'}$ we have $S(n,T)\toop2^{d-r'}$, and by (\ref{2-30})
the variables $V(n,T)$ converge in probability to a limit which
is $[0,\infty)$-valued (actually, it is a.s. positive, but we do not
use this fact here), so that $\De_n\,V(n,T)\toop0$. Since $r'\neq r$,
(\ref{P-70}) readily follows from the definition of $\ca(\al)_T^{n,=r}$.

Propositions \ref{P2} and \ref{P20} are proved analogously, the
alternative-consistency in the latter case following from the second part
of \eqref{2-54}. \qed

\end{document}